\documentclass[12pt,a4paper]{article}

\usepackage[british]{babel}

\usepackage[a4paper,top=2cm,bottom=2cm,left=2.5cm,right=2.5cm,marginparwidth=1.75cm]{geometry}

\usepackage{amsmath}
\usepackage{graphicx}
\usepackage{cite}
\usepackage[colorlinks=true, allcolors=blue]{hyperref}
\usepackage{hyperref}
\usepackage{makecell}
\usepackage[title]{appendix}
\usepackage{mathrsfs}
\usepackage{amsfonts}
\usepackage{booktabs} 
\usepackage{caption}  
\usepackage{threeparttable} 
\usepackage{algorithm}
\usepackage{algorithmicx}
\usepackage{algpseudocode}
\usepackage{listings}
\usepackage{enumitem}
\usepackage{chngcntr}
\usepackage{booktabs}
\usepackage{lipsum}
\usepackage{subcaption}
\usepackage{authblk}
\usepackage[T1]{fontenc}    
\usepackage{csquotes}       
\usepackage{diagbox}
\allowdisplaybreaks[2]

\numberwithin{equation}{section}
\usepackage{setspace}
\usepackage{titlesec}
\titleformat{\section} 
  {\normalfont\Large\bfseries}{\thesection.}{1em}{}
  
\usepackage{lineno} 

\rightlinenumbers 

\linenumbers 

\usepackage{float}   
\usepackage{caption} 


\title{Subsystem Resetting of a Heterogeneous Network of Theta Neurons}

\author[1,2]{Na Zhao}
\author[2]{Carlo R Laing}
\author[1]{Jian Song}
\author[1,*]{Shenquan Liu}
\affil[1]{Department of Mathematics, South China University of Technology, Guangzhou 510641, China}
\affil[2]{School of Mathematical and Computational Sciences, Massey University, Auckland 4442, New Zealand}

\date{}  

\begin{document}
\maketitle

\begin{abstract} 
	Stochastic resetting has shown promise in enhancing the stability and controling of activity in various dynamical systems. In this study, we extend this framework to theta neuron network by exploring the effects of partial resetting, where only a fraction of neurons is intermittently reset. Specifically, we analyze both infinite and finite reset rates, using the averaged firing rate as an indicator of network activity stability. For infinite reset rate, a high proportion of resetting neurons drives the network to stable resting or spiking states. This process collapses the bistable region at the Cusp bifurcation, resulting in smooth and predictable transitions. In contrast, finite resetting introduces stochastic fluctuations, leading to more complex dynamics that occasionally deviate from theoretical predictions. These insights highlight the role of partial resetting in stabilizing neural dynamics and provide a foundation for potential applications in biological systems and neuromorphic computing.
\end{abstract}



\section{Introduction}\label{sec1}

Stochastic resetting, initially introduced within the context of optimization problems and diffusion processes, has garnered significant interest due to its ability to improve convergence in dynamical systems by intermittently returning the system to a specific state. Originally applied to classical diffusion and search algorithms, this mechanism was shown to optimize performance metrics such as mean first-passage time by restarting the system from a predetermined state \cite{evans2011diffusion, Evans2020}. Since then, stochastic resetting has been successfully extended to a wide range of domains \cite{Meylahn2015Large, Masoliver2019, Bodrova2019, Ramoso2020, Mukherjee2018}, including dynamical systems \cite{Ray2021, Singh2022, Jolakoski2023} and networks \cite{Sarkar2022, Bressloff2024chaos, PhysRevE.109.064137}, offering new insights into the behavior of complex systems under resetting protocols.

Most existing studies focus on global resetting, where all components of a system are simultaneously reset. While effective in some scenarios, global resetting becomes impractical for large, interconnected systems due to the difficulty of achieving coordinated resets across all components \cite{Banos2016}. A more scalable alternative is subsystem resetting, where only a fraction of the system is reset at random intervals while the remaining components continue evolving naturally. This selective resetting mechanism, recently applied to the Kuramoto model—a canonical framework for coupled oscillators—has shown that resetting a small subset of oscillators can induce global synchronization, even in systems that do not naturally synchronize \cite{PhysRevE.109.064137}.

Extending subsystem resetting to neural systems offers a promising avenue for regulating complex neural dynamics. Among such systems, theta neuron networks are widely used to model Type-I excitable systems, which exhibit spiking behavior in response to varying excitatory and inhibitory inputs \cite{izhikevich2003simple, Gast2023, luke2013collective}. Unlike the Kuramoto model, which primarily addresses synchronization in coupled oscillators, theta neuron networks provide a biologically grounded framework for studying neural dynamics at the macroscopic level. They have become instrumental in exploring phenomena such as synchronization, oscillations, and firing rate dynamics \cite{Laing2023, Laing2021, Laing2018}. Despite their robust theoretical foundation, the application of stochastic resetting to theta neuron networks remains largely unexplored. While some studies have investigated resetting in neuronal spiking models, most focus on specific input-output characteristics \cite{de2023resetting, chen2023resetting, ermentrout2001phase, Dumont2019, Voloh2015}.

Subsystem resetting provides a scalable alternative by selectively resetting only a fraction of neurons, allowing for targeted control of neural activity. This flexibility enables modulation of the system’s collective behavior, inducing synchronous or asynchronous states depending on parameters such as the fraction of reset neurons and the resetting rate. Such mechanisms hold significant potential for understanding brain dynamics and for designing adaptive artificial neural systems.

Our work introduces a novel application of stochastic resetting within neural networks, examining how partial resetting impacts macroscopic behaviors using the averaged firing rate as an indicator of system stability rather than synchronization. We analyze the regulation of transitions between uniform spiking and resting states across different parameter planes under both finite and infinite reset rates, demonstrating that partial resetting helps stabilize neural activity and reduces large fluctuations in the firing rate. Our results also reveal that finite reset rates introduce stochastic fluctuations, leading to richer and more complex dynamics than traditional models with only global resetting. These findings highlight partial resetting as a versatile mechanism for regulating stability in neural systems. Through comprehensive analytical and numerical analyses, our study provides insights relevant to both biological neuroscience and neuromorphic computing, where controlled resetting mechanisms could enhance the regulation of neural activity.

The remainder of this paper is organized as follows: Section \ref{sec2} presents the theta neuron network with subsystem resetting. Sections \ref{sec3} and \ref{sec4} detail our analytical and numerical findings on firing rate dynamics and bifurcation behavior under different resetting protocols. Finally, in Section \ref{sec5}, we discuss the broader implications of our results and outline potential directions for future research.

\section{The theta neuron network in resetting of subsystem resetting}\label{sec2}

\subsection{The bare theta neuron network}\label{subsec21}
We first give the bare network composed of $N$  theta neurons, all interconnected through pulselike synapses \cite{luke2013complete}. The dynamics of each neuron's phase, $\theta_j$, evolve according to the equation:
\begin{equation}\label{2.1}
	\dot{\theta}_j = (1 - \cos \theta_j) + (1 + \cos \theta_j) [\eta_j + I_{syn}],
\end{equation}
where $j = 1,\dots,N$. The term $I_{syn}$, representing the synaptic input, effectively alters the excitability of the $j$th neuron by modulating its response to external and internal stimuli. The synaptic current is a collective signal generated by the network, with each neuron contributing based on its individual phase. Specifically,
\begin{equation}\label{2.2}
	I_{syn} = \frac{K}{N} \sum_{i}^{N} P_n (\theta_i),
\end{equation}
where $P_n(\theta) = a_n (1 - \cos \theta)^n$, and $n \in \mathbb{N}$ is a parameter controlling the profile of the synaptic pulses. The normalization constant $a_n$ is chosen such that
\begin{equation}\label{2.3}
	\int_{0}^{2\pi} P_n (\theta) \, d\theta = 2\pi,
\end{equation}
ensuring that the total contribution of each neuron integrates correctly over one full cycle.

The global coupling strength $K$ dictates how strongly neurons influence one another, while the parameter $n$ shapes the synaptic pulse, with larger $n$ leading to more peaked synaptic inputs at $\theta=\pi$. This basic setup provides a foundation for understanding how neuron interactions lead to complex network behaviors.

\subsubsection{Heterogeneity}\label{subsubsec211}

For the natural variability (or heterogeneity),
we adopt a Lorentzian distribution to describe the likelihood of any given excitability parameter $\eta$, being present in the network. The probability of selecting an excitability in the range $[\eta, \eta + d\eta]$ is defined by:
\begin{equation}\label{2.4}
	g(\eta) \, d\eta = \frac{\Delta}{\pi \left[ (\eta - \eta_0)^2 + \Delta^2 \right]} \, d\eta.
\end{equation}
where $\eta_0$ is the mean excitability of the population, and $\Delta$ represents the degree of variability, with larger $\Delta$ values indicating greater heterogeneity.

By allowing both positive and negative values of $\eta$, the network consists of a mix of excitable neurons (which require external input to fire) and spontaneously spiking neurons (which fire independently). The balance between these two types is determined by the value of $\eta_0$, which biases the network toward either excitability or spontaneous spiking.
%

\subsubsection{Mean field or order parameter definition}\label{subsubsec212}

To capture the collective behavior of the network, we introduce the mean field or order parameter, a measure often used to describe synchronization in large networks. Following the assumptions outlined in \cite{luke2013collective}: 
\begin{enumerate}
	\item The theta neuron model serves as the canonical form for Type-I excitable systems and their networks.
	\item The network exhibits global (all-to-all) coupling, where each neuron interacts with every other neuron.
	\item  In the large network limit (i.e. as $N\rightarrow\infty$, the thermodynamic limit), we focus on the macroscopic behavior of the network rather than individual neuron dynamics.
\end{enumerate}
The mean field, $z(t)$ , is defined as \cite{kuramoto1975international}:
\begin{equation}\label{2.5}
	z(t) \equiv \frac{1}{N} \sum_{j=1}^{N} e^{\text{i} \theta_j (t)},
\end{equation}
where $z(t)$ provids a measure of its overall rhythm or coherence. This parameter serves as a macroscopic indicator of synchronization within the population, with its magnitude and average phase reflecting the degree of synchrony and the collective firing activity.  

\subsection{The network in presence of subsystem resetting}\label{subsec22}

Building upon the behavior of the bare theta neuron network, we introduce a mechanism of subsystem resetting, inspired by similar approaches used in the Kuramoto model \cite{PhysRevE.109.064137}, where certain neurons are randomly reset to their initial phase. This resetting disrupts the standard evolution of the network, adding a stochastic component that alters the network's long-term dynamics.

We initialize the dynamics from a fully synchronized state, where all neurons have their phase set to $\theta_j(0) = \pi$ for all neurons. At random intervals, the base evolution, governed by Eq. (\ref{2.1}) is interrupted as a selected subset of neurons is reset back to their initial phase $\pi$. This reset process is applied to a predefined subset of neurons, and the same subset undergoes repeated resets across multiple realizations and different reset times. The subset is constructed by uniformly and independently selecting $b$ neurons (where $b<N$) from the total population of $N$ neurons.

We refer to this subset of reset neurons as the \textbf{reset subsystem}, indexed by $j = 1, 2, \dots, b$, while the remaining neurons form the \textbf{non-reset subsystem}, indexed by $j = b+1, b+2, \dots, N$. Resets occur at random times, modeled as a Poisson point process with rate $\lambda$, meaning the time interval $\tau$ between two consecutive resets follows an exponential distribution:
\begin{equation}\label{2.7}
	p(\tau) = \lambda e^{-\lambda \tau}; \quad \lambda \geq 0, \quad \tau \in [0, \infty),
\end{equation}
which reflects the stochastic nature of the resetting events.

Between reset events, during an infinitesimally small interval $dt$ the system evolves according to the intrinsic dynamics described by Eq.(\ref{2.1}) with a probability of $(1-\lambda dt)$. With probability $\lambda dt$, the phases of the $b$ neurons in the reset subsystem are instantaneously reset to $\pi$, while the non-reset neurons continue evolving according to the same dynamics, unaffected by the reset.

In the thermodynamic limit, where $N\rightarrow\infty$, the size of the reset subsystem is described by the fraction $\gamma\equiv b/N$, representing the proportion of neurons subject to resetting. The remaining fraction, $1-\gamma$, corresponds to the size of the non-reset subsystem. Importantly, $\gamma$ can vary between $0$ and $1$, indicating that any proportion of neurons, from none to all, may be reset.

Since the reset and non-reset subsystems follow distinct evolutionary pathways, it becomes essential to track their dynamics separately. To do so, we define the following mean fields to represent the collective behavior of each subsystem:
\begin{subequations}\label{2.8}
	\begin{align}
		z_r(t)
		&\equiv \frac{1}{b} \sum_{j=1}^{b} e^{\text{i} \theta_j (t)}, \label{2.8a} \\
		z_{nr}(t)
		&\equiv \frac{1}{N - b} \sum_{j=b+1}^{N} e^{\text{i} \theta_j (t)}. \label{2.8b}
	\end{align}
\end{subequations}

Here, $z_r(t)$ and $z_{nr}(t)$ represent the complex-valued mean fields of the reset and non-reset subsystems, respectively, capturing the phase coherence and synchronization within each group. By monitoring these quantities, we gain insights into the impact of the resetting mechanism on the overall network dynamics and the interaction between the reset and non-reset neurons, offering a detailed understanding of the system’s emergent behavior.

Next, we explore successively two scenarios for the resetting rate $\lambda$: the infinite and finite cases, following the approach in the Kuramoto model \cite{PhysRevE.109.064137}.

\section{Case 1: Limit $\lambda\rightarrow\infty$}\label{sec3}

In this section,we explore the limiting case where the resetting rate, $\lambda$, approaches infinity. This extreme case provides a useful simplification, allowing us to understand how constant resetting of neurons affects the overall network dynamics.

As $\lambda\rightarrow\infty$, neurons in the reset subsystem remain perpetually pinned at their reset phase, effectively holding their positions at a fixed angle.  Specifically, the phases of these neurons are pinned at $\theta_j(t) = \pi ~ \forall t$, $j = 1,\dots, b$. These neurons, held at their reset phase, cease to evolve dynamically and instead exert a constant influence on the neurons in the non-reset subsystem.

To capture the effect of these pinned neurons on the synaptic input, we substitute their phase condition into Eq. (\ref{2.2}), yielding the following expression for the synaptic current:
\begin{equation}\label{3.1}
	\begin{aligned}
		I_{syn} &= \frac{K}{N} \left[\sum_{i}^{b}a_n (1 - \cos \pi)^n +\sum_{b+1}^{N} a_n (1 - \cos \theta_i)^n\right]\\
		&=Ka_n\left[2^n\gamma +\frac{1}{N} \sum_{b+1}^{N}(1 - \cos \theta_i)^n\right].
	\end{aligned}
\end{equation}
This constant term, $2^n\gamma$, represents the influence of the reset subsystem on the non-reset subsystem through synaptic coupling.

Next, we substitute Eq.(\ref{3.1}) into Eq.(\ref{2.1}), the governing equation for the phase dynamics of the non-reset neurons becomes:
\begin{equation}\label{3.2}
	\begin{aligned}
		\dot{\theta}_j &= (1 - \cos \theta_j) + (1 + \cos \theta_j) \left[\eta_j + Ka_n\left(2^n\gamma +\frac{1}{N}\sum_{b+1}^{N}(1 - \cos \theta_i)^n\right)\right]\\
		&= (1 + \eta_j) - (1 - \eta_j)\cos \theta_j + (1 + \cos \theta_j)Ka_n\Bigg[2^n\gamma +\frac{1}{N}\sum_{b+1}^{N}(1 - \cos \theta_i)^n\Bigg],
	\end{aligned}
\end{equation}
where $j=b+1,\dots,N$. In this form, the infinite resetting rate introduces a constant forcing term, $2^na_nK\gamma$ , that acts on each neuron in the non-reset subsystem. This term essentially represents the aggregate effect of the pinned reset neurons on the rest of the network, modulating the dynamics of the non-reset neurons.

The presence of this constant forcing term allows us to draw an analogy with the classical Kuramoto model, in which a group of oscillators interacts through a mean-field coupling. Following a similar approach to that of the Kuramoto model analysis, we can apply the Ott-Antonsen (OA) ansatz \cite{ott2008low,ott2009long}, a powerful tool for analyzing large networks of coupled oscillators. The OA ansatz reduces the complexity of the system, enabling the derivation of analytical results for the non-reset subsystem, as demonstrated in previous work on similar models \cite{PhysRevE.109.064137}.

\subsection{Mean field reduction of the non-reset subsystem (\ref{3.2})}\label{subsec31}

In this subsection, we explore the reduction of the non-reset subsystem in the thermodynamic limit. By moving to a continuum description, we can express the collective behavior of the neurons using a probability density function, allowing for a mean-field reduction that significantly simplifies the analysis.

\subsubsection{Probability density function}\label{subsubsec311}

In the thermodynamic limit, as the number of neurons  $N\rightarrow\infty$, the behavior of the non-reset subsystem in Eq. (\ref{3.2}) can be described by a probability density function $F(\theta, \eta, t)$, which gives the probability of finding a neuron with phase in the interval $[\theta, \theta + d\theta]$ and excitability in $[\eta, \eta + d\eta]$ at time $t$. $F(\theta, \eta, t)$ must satisfy the normalization condition:
\begin{equation}\label{3.3}
	\int_{-\pi}^{\pi} d\theta \int_{-\infty}^{\infty} d\eta \, F(\theta, \eta, t) = 1,
\end{equation}
holds for all times $t$. We can expand the phase distribution in terms of a Fourier series in $\theta$:
\begin{equation}\label{3.4}
	F(\theta, \eta, t) = g(\eta) \frac{1}{2\pi} \left[ A_0 + \sum_{k=1}^{\infty} \left( A_k(\eta, t) e^{\text{i}k\theta} + \bar{A}_k(\eta, t) e^{-\text{i}k\theta} \right) \right],
\end{equation}
where $A_k$ are the Fourier coefficients that depend on $\eta$ and $t$, and $g(\eta)$ represents the excitability distribution function, assumed to be time-independent:
\begin{equation}\label{3.5}
	g(\eta) = \int_{-\pi}^{\pi} F(\theta, \eta, t) \, d\theta,
\end{equation}
Since the total number of neurons remains conserved, the distribution function $F(\theta, \eta, t)$ satisfies the continuity equation:
\begin{equation}\label{3.6}
	\frac{\partial F(\theta, \eta, t)}{\partial t} + \frac{\partial}{\partial \theta} \left[ F(\theta, \eta, t) \nu(\theta, \eta, t) \right] = 0,
\end{equation}
where $\nu(\theta, \eta, t)$ is the phase velocity of the neurons, which is derived from the continuum version of Eq. (\ref{3.2}).

In the continuum limit, the synaptic current and interactions between neurons can be expressed as integrals over the distribution function $F(\theta, \eta, t)$. The phase velocity $\nu(\theta, \eta, t)$ is given by:
\begin{equation}\label{3.7}
	\begin{aligned}
	\nu(\theta, \eta, t) =& (1+\eta) - (1 - \eta) \cos \theta + K a_n (1 + \cos \theta)\Bigg[ 2^n \gamma + (1  -\gamma)\int_0^{2\pi} d\theta' \int_{-\infty}^{\infty} d\eta' \, \\ 
	& \times F(\theta', \eta', t) (1 - \cos \theta') \Bigg].
	\end{aligned}
\end{equation}
To simplify the higher-order powers of $(1 - \cos \theta')^n$, we apply the binomial theorem and Euler’s representation of the cosine function, yielding:
\begin{equation}\label{3.8}
	(1 - \cos \theta')^n = \sum_{j=0}^{n} \frac{n! (-1)^j}{2^j j! (n-j)!} e^{\text{i}j\theta'} \sum_{m=0}^{j} \frac{j!}{m! (j-m)!} e^{-\text{i} 2m \theta'} 
	= \sum_{j=0}^{n} \sum_{m=0}^{j} P_{jm}e^{\text{i} (j-2m) \theta'}.
\end{equation}
where the coefficients $P_{jm}$ consolidate the factorial terms and are defined as:
\begin{equation}\label{3.9}
	P_{jm} \equiv \frac{n! (-1)^j}{2^j m! (j-m)! (n-j)!}.
\end{equation}
Substituting this result into Eq. (\ref{3.7}), we can rewrite the phase velocity $\nu(\theta, \eta, t)$ as:
\begin{equation}\label{3.10}
		\begin{aligned}
		\nu(\theta, \eta, t) =& (1+\eta) - (1 - \eta) \frac{e^{i\theta} + e^{-i\theta}}{2} + a_nK (1 + \frac{e^{i\theta} + e^{-i\theta}}{2}) \\
		&\times \left[ 2^n \gamma + (1-\gamma) \sum_{j=0}^{n} \sum_{m=0}^{j} P_{jm} \int_0^{2\pi} d\theta' \int_{-\infty}^{\infty} d\eta' \, F(\theta', \eta', t) e^{i(j-2m)\theta'} \right].
	\end{aligned}
\end{equation}

%
In the continuum limit, the discrete mean field of non-reset subsystem in Eq. (\ref{2.8b}) is replaced by an integral over the distribution function:
\begin{equation}\label{3.11}
	z_{nr}(t)\equiv \frac{1}{N - b} \sum_{q=b+1}^{N} e^{\text{i} \theta_q (t)}
	=\int_0^{2\pi} d\theta' \int_{-\infty}^{\infty} d\eta' \, F(\theta', \eta', t) e^{\text{i}\theta'}.
\end{equation}
We also introduce higher-order moments, often referred to as Daido moments \cite{daido1992order,daido1996onset}, defined as:
\begin{equation}\label{3.12}
	z_{nr}^{a}(t) \equiv \int_0^{2\pi} d\theta' \int_{-\infty}^{\infty} d\eta' \, F(\theta', \eta', t) e^{\text{i}a\theta'}.
\end{equation}
Upon inspection of (\ref{3.10}) with $a \equiv j-2m$, and define the continuous influence function $H_n(z_nr)$ and new constant $H$:
\begin{equation}\label{3.14}
	\begin{aligned}
		H_n(z_{nr}) &\equiv a_n \sum_{j=0}^{n} \sum_{m=0}^{j} P_{jm} z_{nr}^{j-2m},\\
		H&\equiv2^na_n.
	\end{aligned}
\end{equation}
Thus, the expression for the phase velocity $\nu$ simplifies to:
\begin{equation}\label{3.15}
	\nu(\theta, \eta, t) = (1+ \eta) - (1 - \eta) \frac{e^{\text{i}\theta} + e^{-i\theta}}{2} + K\left(1 + \frac{e^{\text{i}\theta} + e^{-\text{i}\theta}}{2} \right) [\gamma H + (1-\gamma) H_n(z_{nr})].
\end{equation}

\subsubsection{The Ott-Antonsen reduction method}\label{subsubsec313}

The phase velocity $\nu$ can be expressed in the standard form of a sinusoidally coupled system as:
\begin{equation}\label{3.16}
	\nu = l e^{\text{i}\theta} +h+ \bar{l} e^{-\text{i}\theta},
\end{equation}
where $l$ and $h$ defined as:
\begin{equation}\label{3.17}
	\begin{aligned}
		l &\equiv -\frac{1}{2} \left[ (1-\eta) - K(\gamma H + (1-\gamma) H_n(z_{nr})) \right] = \bar{l},\\
		h &\equiv (1+\eta) + K(\gamma H + (1-\gamma) H_n(z_{nr})).
	\end{aligned}
\end{equation}

Following the OA procedure \cite{ott2008low, ott2009long}, the Fourier coefficients $A_k$ of the oscillator distribution function are expressed in terms of a single complex function $\alpha(\eta, t)$ as:
\begin{equation}\label{3.18}
	A_k = \bar{\alpha}(\eta, t)^k,
\end{equation}
with the condition $|\alpha(\eta, t)| < 1$ at all times $t$. Then, the differential equation for $\alpha$ defined as:
\begin{equation}\label{3.20}
	\begin{aligned}
		\frac{\partial \alpha}{\partial t} &= i (l\alpha^2 + h\alpha + \bar{l}) \\
		&= i [l(\alpha^2+1) + h\alpha].
	\end{aligned} 
\end{equation}

Next, the mean field parameter $z_{nr}(t)$ is expressed as:
\begin{equation}\label{3.22}
	z_{nr}(t) = \int_{-\infty}^{\infty} d\eta' g(\eta') \alpha(\eta', t).
\end{equation}
To proceed, we substitute the specific form of the excitability distribution $g(\eta')$ as a Lorentzian function, given by Eq. (\ref{2.4}):
\begin{equation}\label{3.23}
	z_{nr}(t) = \frac{\Delta}{\pi} \int_{-\infty}^{\infty} \frac{1}{(\eta' - \eta_0)^2 + \Delta^2} \alpha(\eta', t) d\eta'.
\end{equation}
By analytically continuing $\alpha(\eta', t)$ into the upper half of the complex plane and assuming that $|\alpha|$ approaches zero as $\text{Im}(\eta') \rightarrow \infty$, we can apply the residue theorem. This gives:
\begin{equation}\label{3.25}
	z_{nr}(t) =  \alpha(\eta=\eta_0 + i \Delta, t),
\end{equation}
where $\eta = \eta_0 + i \Delta$ is the simple pole within the upper half-plane.

We now return to the expression for the continuous influence function $H_n(z_{nr})$, which incorporates the higher-order Fourier coefficients:
\begin{equation}\label{3.27}
	H_n(z_{nr}) = a_n \left( A_0 + \sum_{q=1}^{n} A_q \left( z_{nr}^q + \bar{z}_{nr}^{q} \right) \right),
\end{equation}
where $A_q = \sum_{j=0}^{n} \sum_{m=0}^{j} \delta_{j-2m,q} P_{jm}$.

Substituting the result for $z_{nr}$ from the residue calculation into the differential equation for $\alpha$ gives us a low-dimensional equation that describes the asymptotic macroscopic behavior of system (\ref{3.2}):
\begin{equation}\label{3.28}
	\frac{d z_{nr}}{d t} = -i \frac{(z_{nr}-1)^2}{2} + [-\Delta + i(\eta_0 +  K(\gamma H + (1-\gamma) H_n(z_{nr})))]\frac{(z_{nr}+1)^2}{2},
\end{equation}

For simplicity and robustness, we set $n = 2$, which has been shown to capture the qualitative dynamics effectively \cite{luke2013collective}. Thus,
\begin{equation}\label{3.29}
	\begin{aligned}
		H_2(z_{nr}) &=\frac{2}{3} \left[\frac{3}{2}-\left( z_{nr} + \bar{z}_{nr} \right)+\frac{\left( z_{nr}^2 + \bar{z}_{nr}^{2} \right)}{4} \right],\\
		H&=2^2(2/3)=8/3.
	\end{aligned}
\end{equation}
Given $\Delta=0.1$, the simple mean field Eq. (\ref{3.28}) is parameterized by three network parameters: the median excitability $\eta_0$, the global coupling strength $K$, and the size of the reset subsystem $\gamma$.

\subsection{Macroscopic dynamics of the non-reset subsystem (\ref{3.2})}\label{subsec32}

In this section, we gain insights into the collective dynamics of the non-reset neurons and focus on understanding the macroscopic behavior of the non-reset subsystem, as governed by Eq. (\ref{3.28}). The macroscopic perspective allows us to connect the behavior of individual neurons to the overall network dynamics, helping us to identify key patterns and states.

\subsubsection{Synchronization}\label{subsubsec320}

We begin by examining whether resetting can synchronize the non-resetting subsystem. To represent the collective dynamics of non-reset subsystem, we employ the order parameter as:
\begin{equation}\label{a}
	z_{nr}(t)=r_{nr}(t) e^{\text{i} \psi_{nr}(t)},
\end{equation}
where real quantity $r_{nr}(t)$  indicates global phase synchrony at time $t$, with $\psi_{nr}(t) \in [0, 2\pi)$ as the average phase. Here, $r_{nr} = 0$ denotes incoherence, $r_{nr} = 1$  perfect synchrony, and $0 < r_{nr} < 1 $ partial synchrony.

Next, applying Eq. (\ref{a}) in Eq. (\ref{3.28}) and comparing the real and imaginary components on both sides, we derive the desired two-dimensional dynamical system:
\begin{equation} \label{a1}
	\begin{aligned}
		\frac{d r_{nr}}{d t} &= \frac{1}{2} (1 - B)(r_{nr}^2 - 1) \sin(\psi_{nr}) - \frac{\Delta}{2} (r_{nr}^2 + 1) \cos(\psi_{nr}) - r_{nr} \Delta, \\
		r_{nr} \frac{d \psi_{nr}}{d t} &= \frac{1}{2} (B - 1)(r_{nr}^2 + 1) \cos(\psi_{nr}) + \frac{\Delta}{2} (1 - r_{nr}^2) \sin(\psi_{nr}) + r_{nr} (1 + B),
	\end{aligned}
\end{equation}
where
\begin{align*} 
	B = \eta_0 + K \left(\frac{8 \gamma}{3} + \frac{2(1 - \gamma)}{3} \left[ \frac{3}{2} - 2 r_{nr} \cos(\psi_{nr}) + \frac{r_{nr}^2 \cos(2 \psi_{nr})}{2} \right]\right).
\end{align*}

\begin{figure}[htp]
	\centering
	\includegraphics[width=0.8\textwidth]{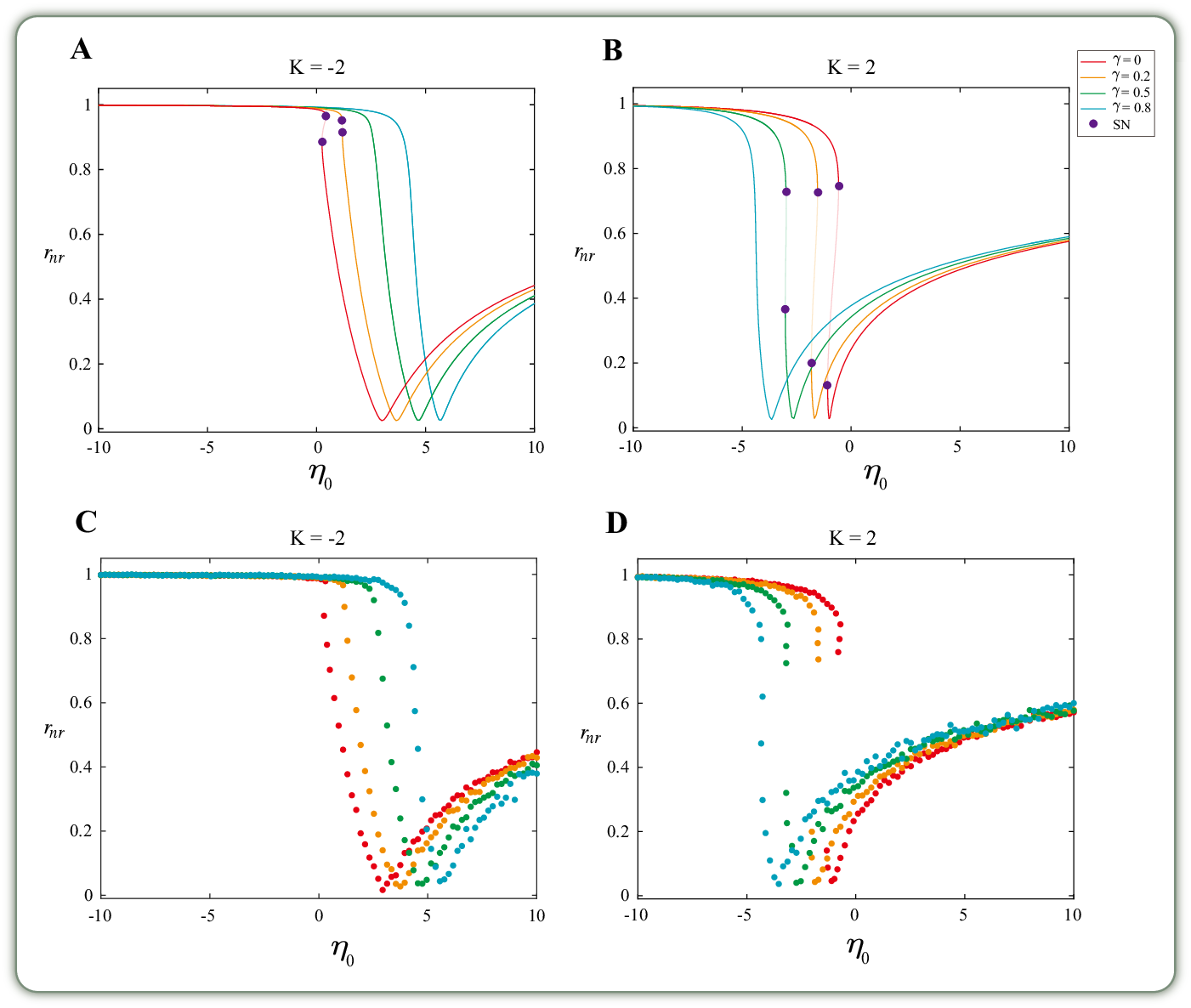}
	\caption{Synchronization results of subsystem resetting with $\lambda\rightarrow\infty, \Delta=0.1$. Plane A-B: Codimension-one bifurcation diagrams with $\eta_0$ as the bifurcation parameter. The SN represent the saddle-node bifurcation. Dark colors represent stable equilibria, while light colors represent unstable equilibria. Plane C-D: Numerical simulation results, the data correspond to a single realization of the dynamics for a system of $N = 10000$ neurons with the integration time step equal to 0.01. See Table \ref{t6} for data on specific bifurcation points.}
	\label{fig:1}       
\end{figure}

Then, we consider both inhibitory ($K=-2.0 < 0$) and excitatory ($ K=2.0 > 0 $) coupling as examples. Figure \ref{fig:1} illustrates the variation of the order parameter $r_{nr}$ with respect to the proportion of resetting neurons ($\gamma$) in the non-resetting subsystem. 

In Fig. \ref{fig:1}, Plane A-B correspond to the theoretical results, using the median excitability $\eta_0$ as the bifurcation parameter, while Plane C-D show the results from numerical simulations of theta neuron network, with simulation details provided in \ref{app1}. Obviously, the theoretical and numerical results align closely. However, in the case of theta neuron network, resetting does not produce a linearly increasing trend in synchronization of the non-resetting subsystem. Initially, global synchrony is maintained for a period due to the initial conditions, but then SN bifurcations trigger the sudden dynamical transitions, with only a brief segment of linear increase appearing at the end. Therefore, resetting does not yield the anticipated outcome of global synchronization for theta neuron network.

Controlling the collective behavior of biological oscillators has increasingly focused on desynchronization strategies, driven by the pathological synchrony observed in disorders such as essential and Parkinsonian tremor. For example, Monga et al. \cite{monga2019phase} developed a control framework to replicate various collective behaviors in biological oscillators. Additionally, Tass \cite{ tass2003model} introduced a desynchronization technique—coordinated resetting of neuronal subpopulations—offering a novel therapeutic approach for conditions characterized by episodic synchronous neuronal oscillations, including Parkinson’s disease, essential tremor, and epilepsy.

Thus, rather than examining how local resetting influences global synchronization in theta neuron networks, as in the Kuramoto model \cite{PhysRevE.109.064137}, we will focus on a newly defined metric, which will serve as the central theme of this study.

\subsubsection{Averaged firing rates}\label{subsubsec321}

It becomes crucial to distinguish between neurons that are quiescent and those that are actively firing. The firing rate of neuron $k$, as introduced in \cite{Omelchenko2022}, is mathematically expressed as:
\begin{equation}\label{2.6}
	f_k = \frac{1}{2\pi} \left\langle \frac{\text{d}\theta_k}{\text{d}t} \right\rangle_T,
\end{equation}
where the time average $\left\langle \cdot \right\rangle_T$  is taken over a long period. 

If all neurons exhibit $f_k=0$ , the network is said to be in a \textbf{uniform rest state}, meaning none of the neurons are firing. Conversely, when all firing rates $f_k$ are positive and equal, the network enters a \textbf{uniform spiking state}, where each neuron fires regularly.

However, identical firing rates do not necessarily imply synchronized firing. Neurons may fire asynchronously or in specific patterns, but after time-averaging, the network exhibits the same mean firing rate across all neurons. This distinction highlights the complexity of neural dynamics, where uniformity in firing rates can mask intricate temporal patterns.

To quantify the firing activity of the non-reset subsystem, we define $W\equiv(1-\bar{z})/(1+\bar{z})$ \cite{Laing2015, Montbrio2015,Omelchenko2019}. The real part of $W$, normalized by $\pi$, gives the firing rate at any given time $t$:
\begin{equation}\label{3.30}
	f(t) = \frac{1}{\pi} \text{Re} \, W = \frac{1 - |z|^2}{\pi |1 + z|^2}.
\end{equation}
If $z(t)$ reaches a steady-state, meaning it no longer depends on time, this formula simplifies to provide the continuum limit analogue of the averaged firing rate, as defined earlier in Eq. (\ref{2.6}). This measure is essential for analyzing whether the neurons in the non-reset subsystem are in a quiescent or actively firing state over time.

\subsubsection{Macroscopic equilibrium states of Averaged firing rates}\label{subsubsec322}

Next, we turn our attention to the long-term behavior of the non-reset subsystem by investigating its equilibrium states. The macroscopic dynamics of the system are governed by the mean-field evolution equation given in (\ref{3.28}). In this setting, the large-scale asymptotic behavior of the non-reset subsystem can be fully captured by a two-dimensional system of ordinary differential equations (ODEs), where the mean-field variable $z_{nr}$ is separated into its real and imaginary components, $z_{nr}\equiv x_{nr} + iy_{nr}$. The resulting dynamical equations are:
\begin{equation}\label{3.31}
		\begin{aligned}
		\dot{x}_{nr} &= (x_{nr} - 1) y_{nr} - \frac{(x_{nr} + 1)^2 - y_{nr}^2}{2} \Delta - (x_{nr} + 1) y_{nr} \left[\eta_0 + K(\frac{8}{3}\gamma + (1-\gamma) \right. \\
		& \left. \times H_2(z_{nr}))\right],\\
		\dot{y}_{nr} &= -\frac{(x_{nr} - 1)^2 - y_{nr}^2}{2} - (x_{nr} + 1) y_{nr} \Delta + \frac{(x_{nr} + 1)^2 - y_{nr}^2}{2} \left[\eta_0 + K(\frac{8}{3}\gamma \right. \\
		& \left. + (1-\gamma) H_2(z_{nr}))\right],
	\end{aligned}
\end{equation}
where $H_2(z_{nr})=\frac{2}{3} \left[\frac{3}{2}-2x_{nr}+\frac{\left( x_{nr}^2 - y_{nr}^{2} \right)}{2} \right]$. 

The averaged firing rate $f_{nr}$ for the non-reset subsystem can be expressed in terms of the mean-field components:
\begin{equation}\label{3.32}
		f_{nr} = \frac{1 - |z_{nr}|^2}{\pi |1 + z_{nr}|^2}
		= \frac{1 - (x_{nr}^2+y_{nr}^2)}{\pi ((1 + x_{nr})^2+ y_{nr}^2)}.\\
	\end{equation}
	
	If $f_{nr}=0$, the corresponding solution of mean field Eq. (\ref{3.28}) means the uniform rest state of non-reset subsystem, by contrast, $f_{nr}>0$, then the corresponding means the uniform spiking state.
	
	Next, we are interested in the equilibrium states of Eq. (\ref{3.31}), and hence, we need to put $\dot{x}_{nr}=0, \dot{y}_{nr}=0$, which in turn implies all solutions to attain a time independent value at long times. Then, the analytical results visualized with the averaged firing rates.
	
	\subsubsection{Varying parameter $\gamma~\&~\eta_0$}\label{subsubsec323}
	
	With the global coupling strength $K$ fixed, we now analyze the effect of resetting on the macroscopic dynamics by varying the median excitability parameter $\eta_0$. This exploration allows us to see how resetting influences the system's transition between rest and spiking states under different excitability conditions. 
	
\begin{figure}[htp]
		\centering
		\includegraphics[width=0.75\textwidth]{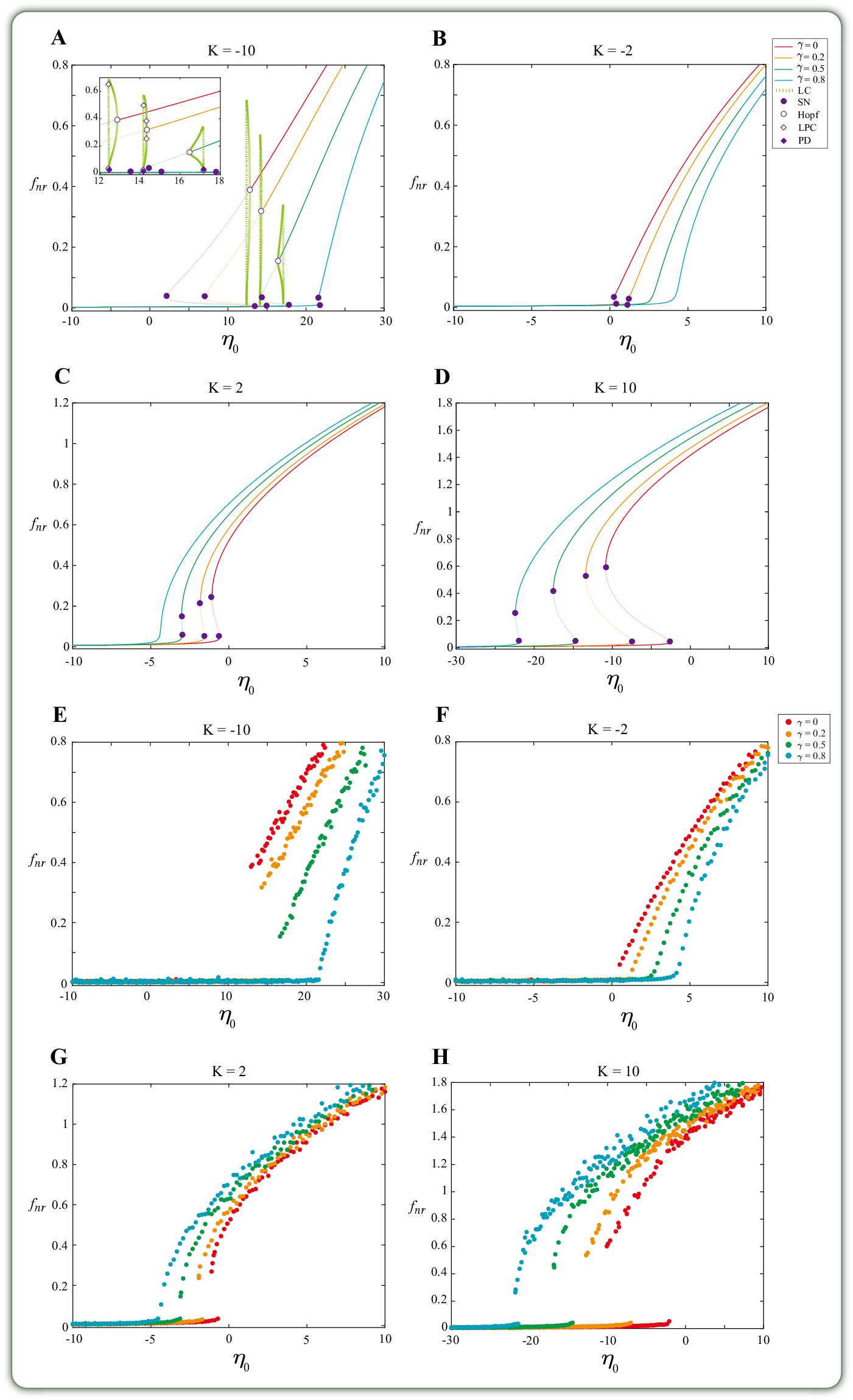}
		\caption{Plane A-D: Codimension-one bifurcation diagrams with $\eta_0$. $\lambda\rightarrow\infty, \Delta=0.1$. LC denotes the limit cycle arising from a Hopf bifurcation. LPC and PD represent the Limit Point bifurcation of cycles and period-doubling bifurcation, respectively. Other legend is the same as Figure \ref{fig:1}. Plane E-H: Numerical simulation results for $\lambda = 50.0$, demonstrating agreement with theory, a single realization of the dynamics for a system of $N = 10000$ neurons with the integration time step equal to 0.01. See Table \ref{t1} and \ref{t5} for data on specific bifurcation points.}
		\label{fig:2}       
	\end{figure}

	In Figure \ref{fig:2} Plane A-D, we examine the bifurcation curves for both inhibitory and excitatory coupling.
	For inhibitory coupling ($K = -10.0, -2.0$), as $\gamma$ increases, the bifurcation curve shifts to the right, indicating that higher values of $\eta_0$ are required to induce spiking. This shift signifies that resetting stabilizes the uniform rest state, expanding the range of $\eta_0$ values over which rest persists. The stronger the inhibitory coupling, the greater the value of $\gamma$ needed to stabilize the system. Notably, as $\gamma$ increases, the averaged firing rate $f_{nr}$ decreases, further reinforcing the rest state.
	
	For $K=-10.0$, the occurrence of Hopf bifurcations leads to the formation of limit cycles, representing periodic solutions within the system (\ref{3.31}). This periodic activity is important as it reflects the collective oscillatory behavior of the theta neuron network. The stability of the limit cycle depends on whether the corresponding Hopf bifurcation is supercritical (stable) or subcritical (unstable) (see Table \ref{t1} for details). The limit cycles ultimately finally
	homoclinic to the saddles at the Limit Point bifurcation of cycles (LPC). The period-doubling bifurcation (PD) may serve as a critical mechanism for the emergence of bursting states \cite{Vo2010, Zhao2024}.
	
	For excitatory coupling ($K = 2.0, 10.0$), however, the situation is somewhat the opposite. Specifically, as $\gamma$ increases, the bifurcation curve shifts to the left, meaning that spiking is induced at lower values of $\eta_0$. This indicates that resetting enhances the uniform spiking state, particularly as $\gamma$ increases, and it increases the averaged firing rate $f_{nr}$. Thus, while resetting suppresses spiking under inhibitory coupling. 
	
	We further validated these theoretical results through numerical simulations of full theta neuron network (\ref{2.1}). Interestingly, even though our theoretical results assume an infinite resetting rate ($\lambda \rightarrow \infty$), the simulations show that a finite rate, as low as $\lambda=50.0$, produces qualitatively similar effects. This demonstrates that numerical simulations with practical values of $\lambda$ are sufficient to capture the resetting dynamics. As shown in the Figure \ref{fig:2} plane E $-$ H, the numerical results align well with the theoretical predictions, confirming that the system behaves consistently across different values of $\gamma$.
	
	\begin{figure}[htp]
		\centering
		\includegraphics[width=0.8\textwidth]{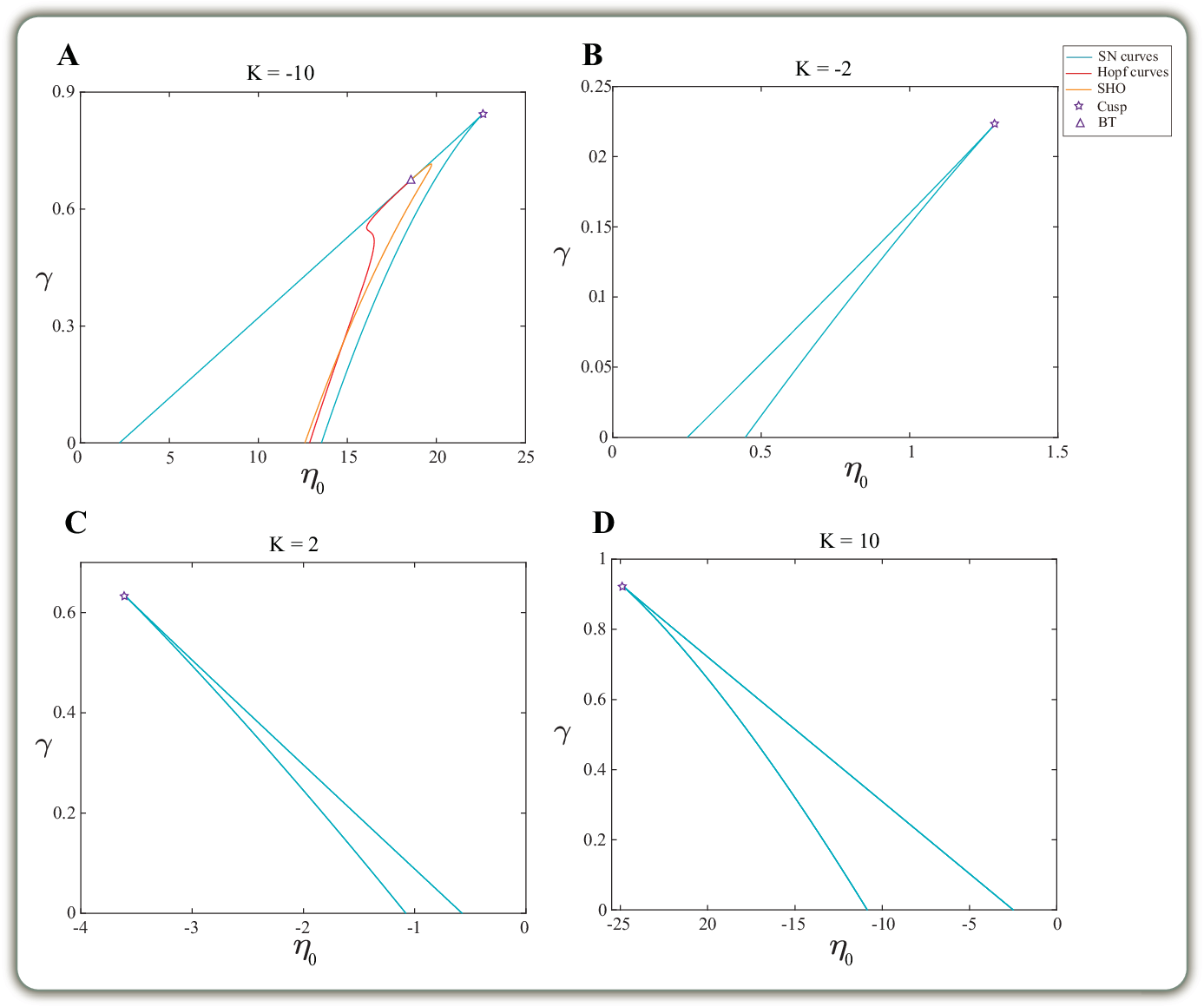}
		\caption{Codimension-two bifurcation in the plane $(\eta_0, \gamma)$.  $\lambda\rightarrow\infty, \Delta=0.1$. The Bogdanov-Takens bifurcations, denoted as BT, signify the tangential intersections between the Hopf curves and the SN curves. The SHO (Saddle Homoclinic Orbit) represents a trajectory where the system's stable and unstable manifolds of a saddle point coincide, emerging after the BT bifurcation. The Cusp bifurcations represent the points where the SN curves intersect. See Table \ref{t2} for data on specific bifurcation points.}
		\label{fig:3}       
	\end{figure}
	Starting from the SN and Hopf bifurcation points in Figure \ref{fig:2}, we further plotted the corresponding SN and Hopf bifurcation curves in the parameter plane $(\eta_0, \gamma)$ (Figure \ref{fig:3}). The region enclosed by the SN bifurcation curves (blue) represents the bistable region in the parameter plane $(\eta_0,\gamma)$. Here, bistability refers to the coexistence of two distinct stable states: uniform rest and uniform spiking. This means that, depending on the initial conditions, system (\ref{3.31}) may converge to different stable states. For instance, as illustrated in Figure \ref{fig:2}, when $\eta_0$ lies between the two SN bifurcations, if the initial state is near the lower branch of the dark stable curve, the solution will converge to this stable equilibrium, indicating uniform rest. Conversely, if the initial state is close to the upper branch of the dark stable curve, the solution will approach this equilibrium, corresponding to uniform spiking. The purple pentagram represents the Cusp bifurcation point, where the two SN bifurcation curves collide and disappear. The red Hopf bifurcation curve is tangent to the SN bifurcation curve at the Bogdanov-Takens (BT) bifurcation, where a saddle homoclinic Orbit (SHO) forms between the Hopf and SN bifurcations.
	
	Figure \ref{fig:3} bifurcation analysis reveals how resetting influences the system's stability under both excitatory and inhibitory coupling conditions. In the case of weaker coupling, only small amounts of resetting are required to eliminate unstable equilibria, whereas stronger coupling demands greater resetting to achieve similar effects.
	
	Thus, when $\eta_0$ is used as the bifurcation parameter, high proportion of resectting neurons can enhance the system's stability and alter the bifurcation phenomena observed under the original parameter control. The larger resetting, the stronger the effect. Interestingly, when inhibitory coupling leads to uniform rest state, resetting expands the range of resting. Conversely, when excitatory coupling leads to uniform spiking, resetting increases the averaged firing rate.
	
	\subsubsection{Varying parameter $\gamma~\&~K$}\label{subsubsec324}
	
	With parameter $\eta_0$ held constant, we now explore how resetting influences the macroscopic dynamics by varying the global coupling strength $K$.
	
\begin{figure}[htp]
		\centering
		\includegraphics[width=0.75\textwidth]{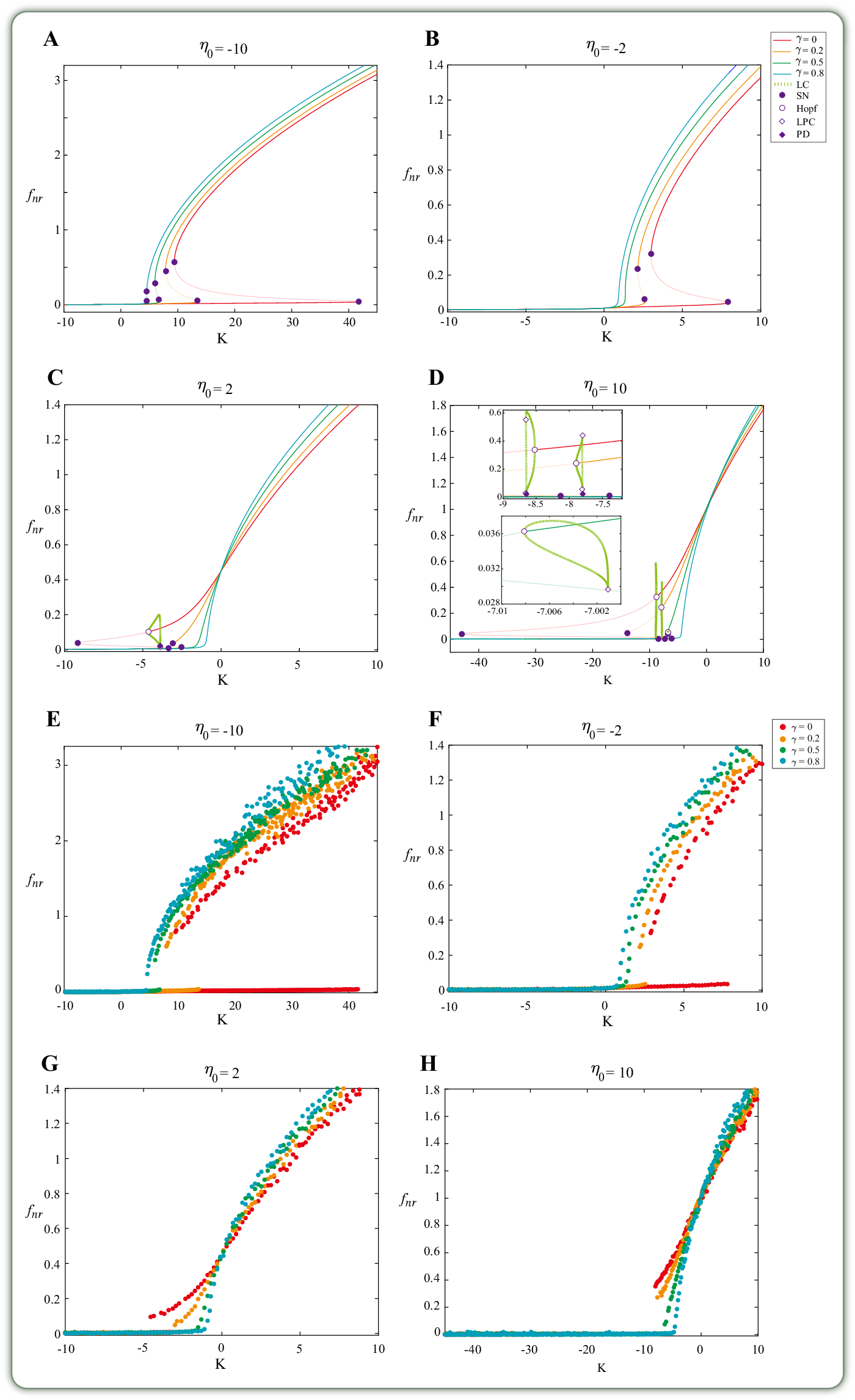}
		\caption{Plane A-D: Codimension-one bifurcation diagrams with $K$. $\lambda\rightarrow\infty, \Delta=0.1$. Legend is the same as Figure \ref{fig:2}. Plane E-H: Numerical simulation results for $\lambda = 50.0$, $N = 10000$ and the integration time step equal to $0.01$. See Table \ref{t1} and \ref{t5} for data on specific bifurcation points.}
		\label{fig:5}       
	\end{figure}

	We treat $K$ as a bifurcation parameter, Figure \ref{fig:5} highlights how the bifurcation structure changes with different values of $\gamma$. For $\eta_0 = -10.0, -2.0<0$, increasing $\gamma$ leads to a gradual leftward shift in the bifurcation curve. This shift means that the transition to spiking occurs at lower $K$ values as the resetting level increases, expanding the range of stable spiking states. Additionally, as $\gamma$ continues to grow, the bifurcations eventually disappear, leaving only stable spiking behavior. 
	
	For $\eta_0 = 2.0, 10.0>0$, we observe a different pattern. Similar limit cycles arising from Hopf bifurcations will not be further elaborated here. Resetting pushes the bifurcation curve to the right, meaning that higher values of $K$ are needed to induce spiking. As a result, resetting suppresses the spiking state for inhibitory coupling ($K<0$), while in the excitatory regime ($K>0$), resetting enhances spiking by increasing the averaged firing rate $f_{nr}$. This nuanced behavior shows how resetting’s effect on network dynamics varies depending on the coupling regime. 
	
	To confirm the accuracy of our theoretical predictions, we validated these results with numerical simulations of full theta neuron network (\ref{2.1}). As illustrated in Figure \ref{fig:5} Plane E-H, which closely match the theoretical results. This highlights the robustness, showing that finite $\lambda$ values still capture the essential dynamics predicted by theory. 
	
\begin{figure}[htp]
		\centering
		\includegraphics[width=0.8\textwidth]{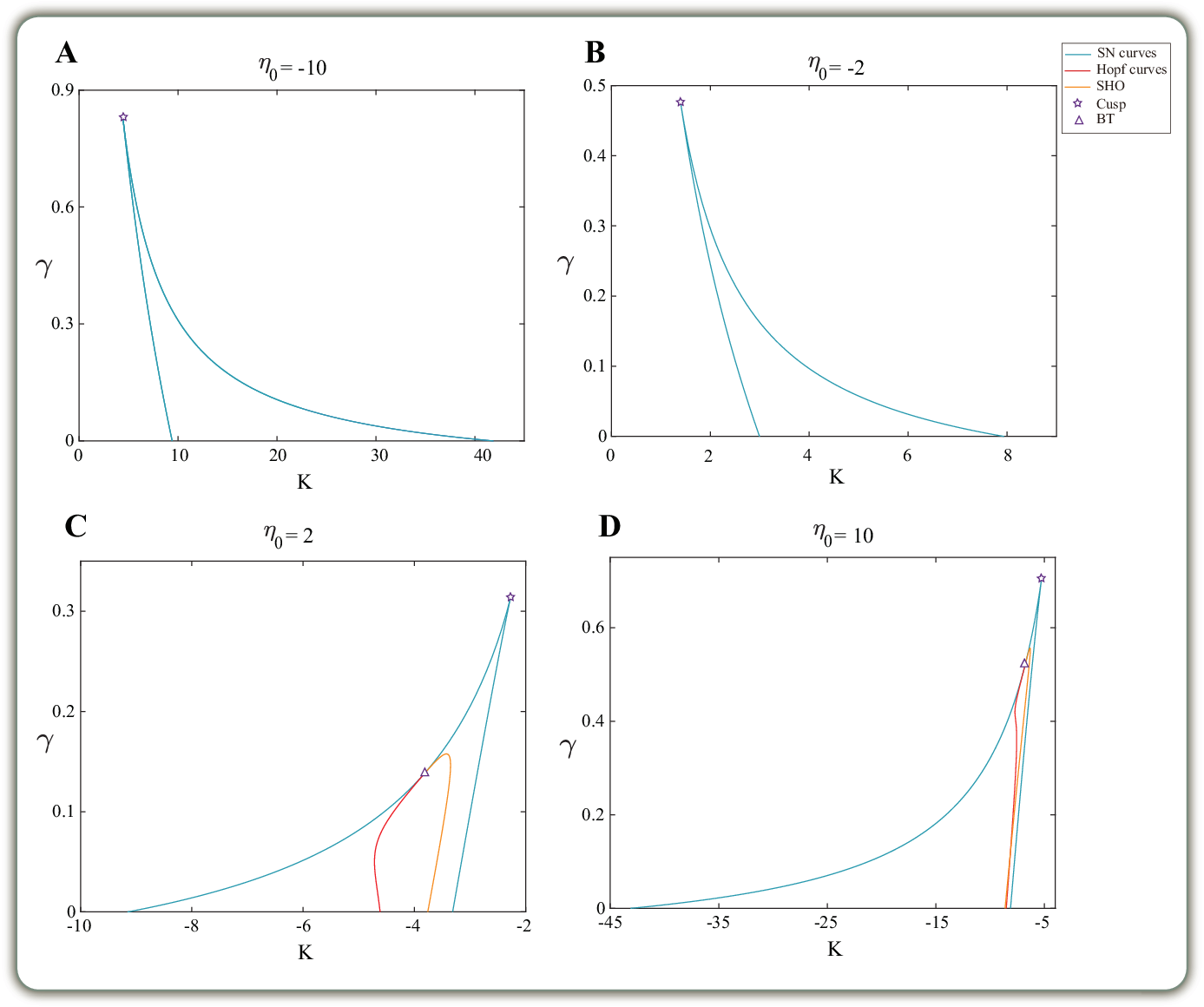}
		\caption{Codimension-two bifurcation in the plane $(K, \gamma)$. $\lambda\rightarrow\infty, \Delta=0.1$. Legend is the same as Figure \ref{fig:3}. See Table \ref{t2} for data on specific bifurcation points.}
		\label{fig:6}       
	\end{figure}

	Lastly, we also plotted the corresponding SN and Hopf bifurcation curves in the plane $(K, \gamma)$, as shown in Figure \ref{fig:6}. Similar bifurcation curves also define corresponding bistable regions and associated codimension-two bifurcation points. For further details, refer to Table \ref{t2}. These critical transitions mark significant changes in the system’s stability, with resetting playing a key role in eliminating unstable equilibria and enhancing stable spiking behavior.

	In summary, when $K$ is the bifurcation parameter, resetting begins to influence the system in the excitatory coupling region ($K > 0$) once the median excitability $\eta_0$ is negative. This is due to the intrinsic properties of the theta neuron model, which promotes firing in this regime. However, when $\eta_0 > 0$, resetting has the opposite effect, where it decreases the averaged firing rate in the inhibitory region but promotes spiking in the excitatory region. Thus, resetting proves to be a versatile tool, capable of both suppressing and enhancing spiking behavior depending on the context of the coupling and excitability parameters.
	
	\section{Case 2: Finite $\lambda$}\label{sec4}
	
	We now consider the more practical scenario of finite resetting rate $\lambda$. Unlike the case with infinite $\lambda$, where neurons in the reset subsystem are effectively fixed at $\theta=\pi$, a finite $\lambda$ introduces intermittent resets, leading to a more complex interplay between the reset and non-reset neurons. This intermittent resetting not only influences the behavior of individual neurons but also alters the collective dynamics of the network as a whole.
	
	In the continuum limit ($N\rightarrow\infty$), the Ott-Antonsen (OA) ansatz provides significant simplifications for analyzing the system under bare evolution \cite{ott2008low}. However, when resetting occurs at a finite rate, the overall probability distribution of the neuronal phases no longer satisfies the OA ansatz, adding further complexity to the analysis.
	
	\subsection{Analysis in absence of resetting}\label{subsec41}
	
	To establish a basis for understanding the effects of finite resetting, we begin by examining the case of bare evolution, in which the network operates without resetting. Under these conditions, we consider the system as comprising two subsystems: reset and non-reset neurons. If the initial phase distribution of each subsystem satisfies the OA ansatz, we can apply the OA ansatz separately to each subsystem to capture the time evolution of their respective mean fields. The dynamics of the subsystems are coupled, with each subsystem's evolution influenced by the other \cite{PhysRevE.109.064137}.
	
	In the absence of resetting, the modified network model (\ref{2.1}) for the reset and non-reset subsystems can be written as:
	\begin{equation}\label{4.1}
		\begin{aligned}
			\frac{d\theta_i^{(r)}}{dt} &= 1 + \eta_i^{(r)} - \left( 1 - \eta_i^{(r)} \right)\cos \theta_i^{(r)} + Ka_n \left( 1 + \cos \theta_i^{(r)} \right) \\
			&\quad \left[ \gamma \sum_{j=1}^{b} \left( 1 - \cos \theta_j^{(r)} \right)^n + \frac{1}{N} \sum_{j=b+1}^{N} \left( 1 - \cos \theta_j^{(nr)} \right)^n \right], \\
			\frac{d\theta_i^{(nr)}}{dt} &= 1 + \eta_i^{(nr)} - \left( 1 - \eta_i^{(nr)} \right)\cos \theta_i^{(nr)} + Ka_n \left( 1 + \cos \theta_i^{(nr)} \right) \\
			&\quad \left[ \gamma \sum_{j=1}^{b} \left( 1 - \cos \theta_j^{(r)} \right)^n + \frac{1}{N} \sum_{j=b+1}^{N} \left( 1 - \cos \theta_j^{(nr)} \right)^n \right].
		\end{aligned}
	\end{equation}
	Here, the superscripts $(r)$ and $(nr)$ indicate neurons belonging to the reset and non-reset subsystems, respectively. The two subsystems interact through synaptic inputs that are influenced by the collective activity of both subsystems. We assume that the sharpness parameter $n$ is identical for both subsystems. 
	
	\subsubsection{Mean field reduction of the two subsystems (\ref{4.1})}\label{subsubsec411}
	
	Using the approach described in Section \ref{subsec31}, we derive a set of low-dimensional dynamical equations to describe the asymptotic macroscopic states of each subsystem. We assume that both subsystems are in the thermodynamic limit, $b\rightarrow\infty$ and $N\rightarrow\infty$, which allows us to represent the probability distributions for each subsystem using continuous functions: $F_{r}(\theta_{r}, \eta_{r}, t)$ and $F_{nr}(\theta_{nr}, \eta_{nr}, t)$. 
	
	Since the reset and non-reset subsystems are considered independent, their respective probability distributions must obey their own continuity equations:
	\begin{equation}\label{4.2}
		\begin{aligned}
			&\frac{\partial F_{r}}{\partial t} + \frac{\partial}{\partial \theta_{r}} \left( F_{r} \nu_{r} \right) = 0, \\
			&\frac{\partial F_{nr}}{\partial t} + \frac{\partial}{\partial \theta_{nr}} \left( F_{nr} \nu_{nr} \right) = 0,
		\end{aligned}
	\end{equation}
	where the phase velocities $\nu_{r}$ and $\nu_{nr}$ for each subsystem are given by:
	\begin{equation}\label{4.4}
			\begin{aligned}
			\nu_{r} &= \frac{e^{i\theta_{r}} + e^{-i\theta_{r}}}{2} \left[ 1 + \eta_{r} + K\gamma H_n(z_{r}) + K(1-\gamma) H_n(z_{nr}) \right] \\
			&\quad + \left[ -(1 - \eta_{r}) + K\gamma H_n(z_{r}) + K(1-\gamma) H_n(z_{nr}) \right] , \\
			\nu_{nr} &= \frac{e^{i\theta_{nr}} + e^{-i\theta_{nr}}}{2} \left[ 1 + \eta_{nr} + K\gamma H_n(z_{r}) + K(1-\gamma) H_n(z_{nr}) \right]  \\
			&\quad + \left[ -(1 - \eta_{nr}) + K\gamma H_n(z_{r}) + K(1-\gamma) H_n(z_{nr}) \right] .
		\end{aligned}
	\end{equation}
	where $H_n$ is defined as (\ref{3.14}).
	
	To simplify further, the phase velocities can be written in a sinusoidally coupled form: $\nu_{p} = l_{p}e^{\text{i}\theta}+ h_{p} + \bar{l}_{p} e^{-\theta}$, where $p=r,nr$ represents the reset and non-reset subsystems, respectively. The coefficients $l_{p}$ and $h_{p}$ defined as:
	\begin{align*}\label{4.5}
		l_p &\equiv -\frac{1}{2} \left[ (1 - \eta_p) - \left( K\gamma H_n(z_{r}) + K(1-\gamma) H_n(z_{nr}) \right) \right] = \bar{l}_p, \\
		h_p &\equiv (1 + \eta_p) + \left( K\gamma H_n(z_{r}) + K(1-\gamma) H_n(z_{nr}) \right),
	\end{align*}
	Using the sinusoidally coupled form of $\nu_{p}$ and following the Ott-Antonsen reduction method, we derive the equations governing the evolution of the Ott-Antonsen parameters $\alpha_{p}$ for each subsystem \cite{Martens2009}:
	\begin{equation}\label{4.6}
		\frac{d\alpha_p}{dt}
		= i \left[ l_p (1 + \alpha_p^2) + h_p \alpha_p \right],
	\end{equation}
	where $|\alpha_p|< 1$ for all time $t$. 
	
	Next, we define the excitability distribution for both subsystems to be the same:
	\begin{equation}\label{4.7}
		g_{r}(\eta) = g_{nr}(\eta) = g(\eta),
	\end{equation}  
	with the median excitability $\eta_0$ and half-width $\Delta$. Using contour integration, we find the relationship between the mean field $z_p(t)$ and the Ott-Antonsen parameter:
	\begin{equation}\label{4.9}
		z_p = \alpha_p \left( \eta_p = \eta_{0} + i \Delta \right),
	\end{equation}
	
	The resulting equations for the time evolution of the mean fields of each subsystem are (also setting $n=2$):
	\begin{subequations}\label{4.12}
		\begin{align}
			\frac{dz_{r}}{dt} &= -i \frac{(z_{r} - 1)^2}{2} + [ -\Delta + i ( \eta_0 +  K\gamma H_2(z_{r}) + K(1-\gamma) H_2(z_{nr}) )]  \frac{(z_{r} + 1)^2}{2}, \label{4.12a} \\
			\frac{dz_{nr}}{dt} &= -i \frac{(z_{nr} - 1)^2}{2} + [ -\Delta + i ( \eta_0 +  K\gamma H_2(z_{r}) + K(1-\gamma)H_2(z_{nr}) )]   \frac{(z_{nr} + 1)^2}{2}, \label{4.12b}
		\end{align}
	\end{subequations}
	where
	\begin{align*}\label{4.11}	
		H_2(z_{r}) &=\frac{2}{3} \left[\frac{3}{2}-\left( z_{r} + \bar{z}_{r} \right)+\frac{\left( z_{r}^2 + \bar{z}_{r}^{2} \right)}{4} \right],\\
		H_2(z_{nr}) &=\frac{2}{3} \left[\frac{3}{2}-\left( z_{nr} + \bar{z}_{nr} \right)+\frac{\left( z_{nr}^2 + \bar{z}_{nr}^{2} \right)}{4} \right].
	\end{align*}
	
	\subsection{Analysis in resetting}\label{subsec42}
	
	We now incorporate resetting into the system dynamics, as described in Section \ref{subsec22}. With finite resetting rate $\lambda$, the mean fields, $z_{r}(t)$ and $z_{nr}(t)$ become random variables, as opposed to the deterministic behavior observed for infinite $\lambda$. In this scenario, it is more appropriate to consider the realization-averaged quantities $\tilde{z}_{r}(t)$ and $\tilde{z}_{nr}(t)$.
	
	At any given time $t$, the mean fields $z_{r}(t)$ and $z_{nr}(t)$ evolve according to Eqs. (\ref{4.12a}) and (\ref{4.12b}). In the interval $[t,t+dt]$, $z_{r}$ evolves either according to its original dynamics (with probability $1 - \lambda dt$) or is reset to $z_{r}(t)=\frac{1}{b} \sum_{j=1}^{b} e^{\text{i} \pi}= -1$ (with probability $\lambda dt$). The average change in the mean fields over time is given by:
	\begin{equation}\label{4.13}
		\begin{aligned}
			d\tilde{z}_r &= (1 - \lambda dt) dz_{r} + \lambda dt(-1 - z_{r}), \\
			d\tilde{z}_{nr} &= dz_{nr}.
		\end{aligned}
	\end{equation}
	where the tilde denotes average over all those realizations. Then, we consider the average of these random variables $z_{r}$ and $z_{nr}$. This gives
\begin{subequations}\label{4.15}
	\begin{align}
		\frac{d\bar{z}_{r}}{dt} =& -i \frac{(\overline{z_{r}^2} - 2\bar{z}_{r} + 1)}{2} + \left[ -\Delta + i \left( \eta_0 +  K\gamma \frac{2}{3} \left( \frac{3}{2} - ( \bar{z}_{r} + \bar{z}_{r}^{*}) + \frac{(\overline{z_{r}^2} + \overline{z_{r}^{*2}})}{4} \right) \right.\right. \nonumber \\
		&\left.\left.+ K(1-\gamma) \frac{2}{3} \left( \frac{3}{2} - ( \bar{z}_{nr} + \bar{z}_{nr}^{*}) + \frac{(\overline{z_{nr}^2} + \overline{z_{nr}^{*2}})}{4} \right) \right) \right] \frac{(\overline{z_{r}^2} + 2\bar{z}_{r} + 1)}{2} \nonumber \\
		&- \lambda(1 + \bar{z}_{r}), \label{4.15a} \\
		\frac{d\bar{z}_{nr}}{dt} =& -i \frac{(\overline{z_{nr}^2} - 2\bar{z}_{nr} + 1)}{2} + \left[ -\Delta + i \left( \eta_0 +  K\gamma \frac{2}{3} \left( \frac{3}{2} - ( \bar{z}_{r} + \bar{z}_{r}^{*}) + \frac{(\overline{z_{r}^2} + \overline{z_{r}^{*2}})}{4} \right) \right.\right. \nonumber \\
		&\left.\left.+ K(1-\gamma) \frac{2}{3} \left( \frac{3}{2} - ( \bar{z}_{nr} + \bar{z}_{nr}^{*}) + \frac{(\overline{z_{nr}^2} + \overline{z_{nr}^{*2}})}{4} \right) \right) \right] \frac{(\overline{z_{nr}^2} + 2\bar{z}_{nr} + 1)}{2}. \label{4.15b}
	\end{align}
\end{subequations}
	
	For large $\lambda$, assuming the distributions of $z_{r}$ and $z_{nr}$are sharply peaked, we can approximate:
	\begin{equation}\label{4.16}
		\begin{aligned}
			\widetilde{z_{r}^2} \approx \tilde{z}_{r}^2,\qquad \widetilde{\bar{z}_{r}^{2}} \approx \tilde{\bar{z}}_{r}^{2},\qquad \widetilde{z_{nr}^2} \approx \tilde{z}_{nr}^2,\qquad \widetilde{\bar{z}_{nr}^{2}} \approx \tilde{\bar{z}}_{nr}^{2},
		\end{aligned}
	\end{equation}
	leading to:
	\begin{subequations}\label{4.17}
		\begin{align}
			\frac{d\bar{z}_{r}}{dt} =& -i \frac{(\bar{z}_{r}^2 - 2\bar{z}_{r} + 1)}{2} + \left[ -\Delta + i \left( \eta_0 +  K\gamma \frac{2}{3} \left( \frac{3}{2} - ( \bar{z}_{r} + \bar{z}_{r}^{*}) + \frac{(\bar{z}_{r}^2 + \bar{z}_{r}^{*2})}{4} \right) \right.\right. \nonumber \\
			&\left.\left.+ K(1-\gamma) \frac{2}{3} \left( \frac{3}{2} - ( \bar{z}_{nr} + \bar{z}_{nr}^{*}) + \frac{(\bar{z}_{nr}^2 + \bar{z}_{nr}^{*2})}{4} \right) \right) \right] \frac{(\bar{z}_{r}^2 + 2\bar{z}_{r} + 1)}{2} \nonumber \\
			&- \lambda(1 + \bar{z}_{r}), \label{4.17a} \\
			\frac{d\bar{z}_{nr}}{dt} =& -i \frac{(\bar{z}_{nr}^2 - 2\bar{z}_{nr} + 1)}{2} + \left[ -\Delta + i \left( \eta_0 +  K\gamma \frac{2}{3} \left( \frac{3}{2} - ( \bar{z}_{r} + \bar{z}_{r}^{*}) + \frac{(\bar{z}_{r}^2 + \bar{z}_{r}^{*2})}{4} \right) \right.\right. \nonumber \\
			&\left.\left.+ K(1-\gamma) \frac{2}{3} \left( \frac{3}{2} - ( \bar{z}_{nr} + \bar{z}_{nr}^{*}) + \frac{(\bar{z}_{nr}^2 + \bar{z}_{nr}^{*2})}{4} \right) \right) \right] \frac{(\bar{z}_{nr}^2 + 2\bar{z}_{nr} + 1)}{2}. \label{4.17b}
		\end{align}
	\end{subequations}
	
	\subsection{Macroscopic dynamics of the two subsystems (\ref{4.1})}\label{subsec43}
	
	To facilitate a more intuitive understanding of the dynamics, we begin by explicitly separating the complex variable $\tilde{z}_{p}\equiv \tilde{x}_{p} + i\tilde{y}_{p}$ into its real and imaginary parts for both the reset ($p=r$) and non-reset ($p=nr$) subsystems. This leads to the following system of ordinary differential equations (ODEs) that govern the macroscopic behavior of the two subsystems (\ref{4.17}):
\begin{equation}\label{4.18}
	\begin{aligned}
		\dot{\bar{x}}_{r} &= ( \bar{x}_{r} - 1) \bar{y}_{r} - \frac{(\bar{x}_{r} + 1)^2 - \bar{y}_{r}^2}{2} \Delta - (\bar{x}_{r} + 1) \bar{y}_{r}  \left[\eta_0 +  K\gamma \frac{2}{3} \left( \frac{3}{2} - 2\bar{x}_{r} \right. \right. \\
		&\quad \left.\left. + \frac{(\bar{x}_{r}^2 - \bar{y}_{r}^2)}{2} \right) + K(1-\gamma) \frac{2}{3} \left( \frac{3}{2} - 2\bar{x}_{nr}+ \frac{(\bar{x}_{nr}^2 - \bar{y}_{nr}^2)}{2} \right) \right]-\lambda(1+\bar{x}_{r}), \\
		\dot{\bar{y}}_{r} &= -\frac{(\bar{x}_{r} - 1)^2 - \bar{y}_{r}^2}{2} - (\bar{x}_{r} + 1) \bar{y}_{r} \Delta + \frac{(\bar{x}_{r} + 1)^2 - \bar{y}_{r}^2}{2}  \left[\eta_0 +  K\gamma \frac{2}{3} \left( \frac{3}{2} - 2\bar{x}_{r} \right. \right. \\
		&\quad \left. \left. + \frac{(\bar{x}_{r}^2 - \bar{y}_{r}^2)}{2} \right) + K(1-\gamma) \frac{2}{3} \left( \frac{3}{2} - 2\bar{x}_{nr}+ \frac{(\bar{x}_{nr}^2 - \bar{y}_{nr}^2)}{2} \right) \right]-\lambda\bar{y}_{r}, \\
		\dot{\bar{x}}_{nr} &= (\bar{x}_{nr} - 1) \bar{y}_{nr} - \frac{(\bar{x}_{nr} + 1)^2 - \bar{y}_{nr}^2}{2} \Delta - (\bar{x}_{nr} + 1) \bar{y}_{nr}  \left[\eta_0 +  K\gamma \frac{2}{3} \left( \frac{3}{2} - 2\bar{x}_{r} \right. \right. \\
		&\quad \left. \left. +\frac{(\bar{x}_{r}^2 - \bar{y}_{r}^2)}{2} \right) + K(1-\gamma) \frac{2}{3} \left( \frac{3}{2} - 2\bar{x}_{nr}+ \frac{(\bar{x}_{nr}^2 - \bar{y}_{nr}^2)}{2} \right) \right], \\
		\dot{\bar{y}}_{nr} &= -\frac{(\bar{x}_{nr} - 1)^2 - \bar{y}_{nr}^2}{2} - (\bar{x}_{nr} + 1) \bar{y}_{nr} \Delta + \frac{(\bar{x}_{nr} + 1)^2 - \bar{y}_{nr}^2}{2}  \left[\eta_0 +  K\gamma \frac{2}{3} \left( \frac{3}{2} \right. \right. \\
		&\quad \left. \left. - 2\bar{x}_{r}  + \frac{(\bar{x}_{r}^2 - \bar{y}_{r}^2)}{2} \right) + K(1-\gamma) \frac{2}{3} \left( \frac{3}{2} - 2\bar{x}_{nr}+ \frac{(\bar{x}_{nr}^2 - \bar{y}_{nr}^2)}{2} \right) \right].
	\end{aligned}
\end{equation}
	
	The system described above captures the large-scale dynamics of the reset ($r$) and non-reset ($nr$) subsystems, with each subsystem’s behavior dictated by a four-dimensional system. In contrast to the case where the reset rate $\lambda \rightarrow \infty$, we now treat $\lambda$ as a tunable parameter, allowing us to investigate the system’s behavior under finite reset rates. In this analysis, we focus primarily on the system’s equilibrium states and the associated firing rates, specifically the averaged firing rate $f_{nr}$ of the non-reset subsystem (\ref{3.32}). 
	
	\subsubsection{Varying parameter $\gamma~\&~\eta_0$}\label{subsubsec431}
	
\begin{figure}[htp]
		\centering
		\includegraphics[width=0.8\textwidth]{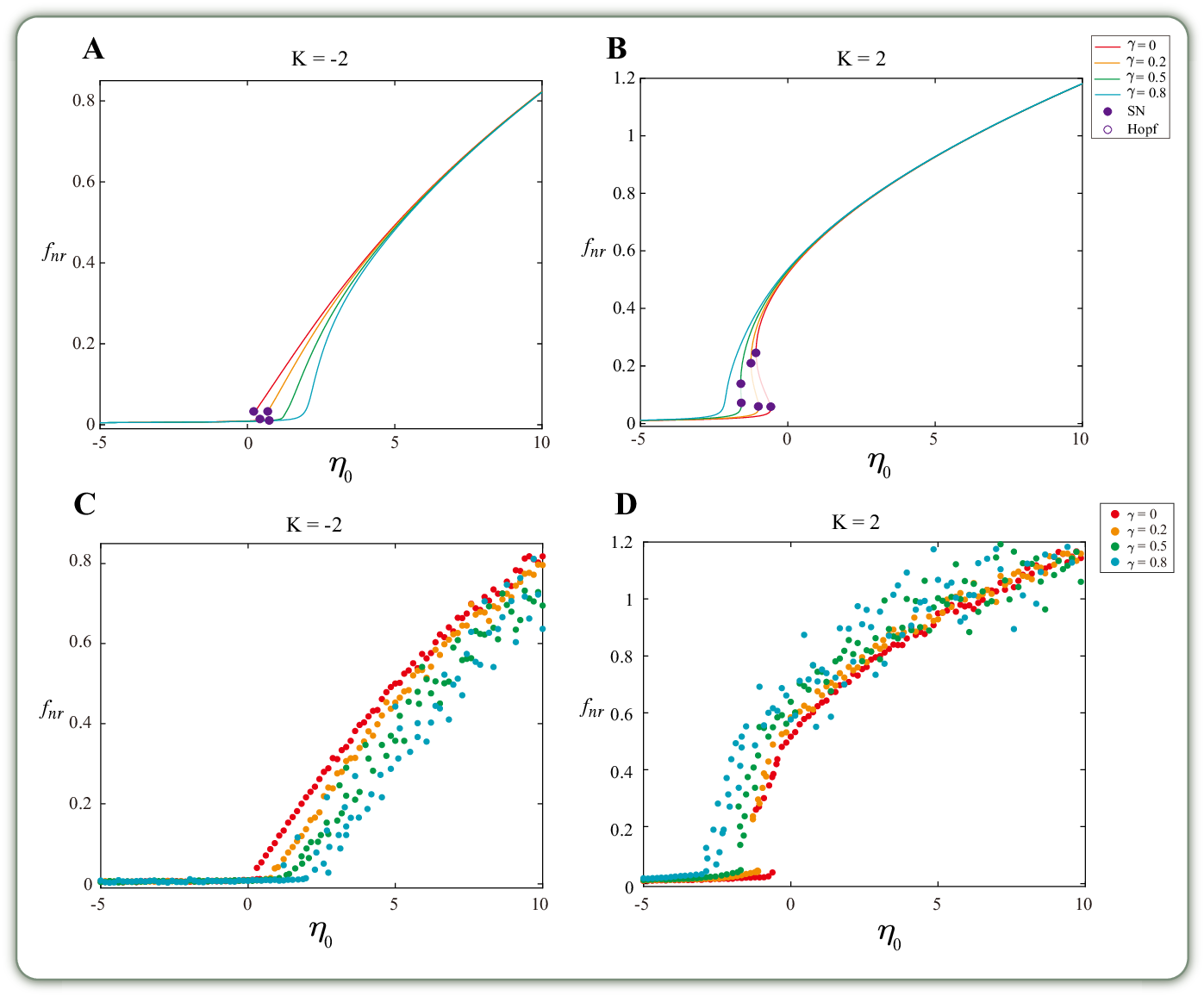}
		\caption{Plane A-B: Codimension-one bifurcation diagrams with $\eta_0$. $\lambda=1.0, \Delta=0.1$. Legend is the same as Figure \ref{fig:2}. Plane C-D: Numerical simulation results. The data correspond to the equilibrium-state value averaged over $100$ realizations of the dynamics for a system of $N = 10000$ neurons with the integration time step equal to $0.05$. See Table \ref{t3} for data on specific bifurcation points.
		}
		\label{fig:7}       
	\end{figure}

	We begin by setting the reset rate $\lambda=1.0$, representing a finite and relatively small value. This allows us to separately explore the effects of excitatory and inhibitory coupling on the system, with a focus on how the parameter $\gamma$ influences the bifurcation behavior with respect to $\eta_0$. The corresponding results are shown in Figure \ref{fig:7}. Despite the small reset rate, $\gamma$ continues to play a significant role in enhancing the stability of the equilibrium. Notably, the transition from uniform rest to uniform spiking follows a trend similar to the case of $\lambda \rightarrow \infty$ (as shown in Figure \ref{fig:2} B and C). Specifically, for inhibitory coupling, the critical transition value increases as $\gamma$ rises, while for excitatory coupling, the critical value decreases. However, these changes are noticeably smaller in magnitude compared to the infinite reset rate scenario.
	
	Interestingly, as the averaged firing rate $f_{nr}$ approaches approximately $0.6$, the equilibrium curve begins to align with the base evolution ($\gamma = 0$), deviating from the behavior observed under infinite reset conditions. This alignment demonstrates that the final trend of the averaged firing rate remains consistent with the non-reset subsystem, even when a large proportion of neurons (up to $80\%$) are being reset under a finite, relatively small reset rate.

	Further supporting this analysis, the numerical results for complete theta neuron network (\ref{2.1}), shown in Figure \ref{fig:7} (Panels C and D), align closely with the theoretical predictions in the uniform rest region, even at the critical transition points. For the uniform spiking state, the numerical results generally agree with theoretical predictions when the proportion of reset neurons is small ($20\%,$ orange. However, when the proportion of reset neurons exceeds half ($50\%-80\%$, green and blue), the numerical simulations exhibit significant jumps with no discernible pattern. This leads to deviations from the theoretical curve, suggesting that the approximation used (Eq. \ref{4.16}) may not be applicable for lower reset rates.
	
\begin{figure}[!htp]
		\centering
		\includegraphics[width=0.8\textwidth]{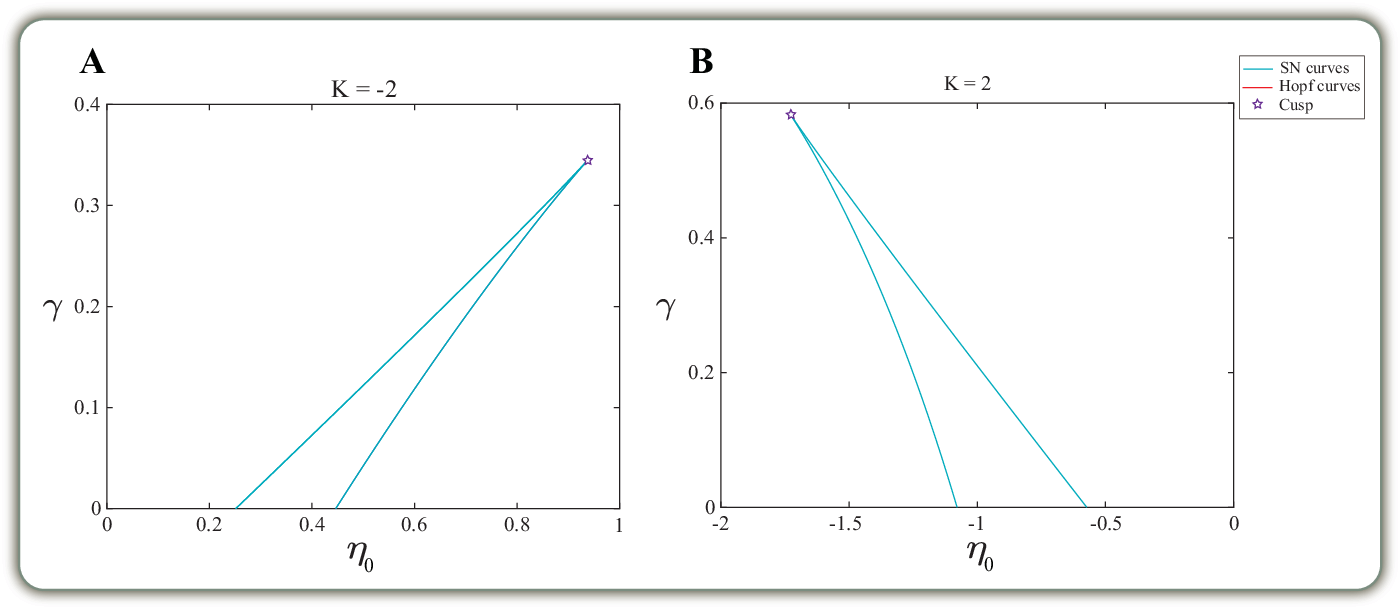}
		\caption{Codimension-two bifurcation in the plane $(\eta_0, \gamma)$. $\lambda=1.0, \Delta=0.1$. Legend is the same as Figure \ref{fig:3}. See Table \ref{t4} for data on specific bifurcation points.}
		\label{fig:8}       
	\end{figure}

	Next, we plot the bifurcation curves in the codimension-two plane $(\eta_0, \gamma)$ also defined bistable regions, as shown in Figure \ref{fig:8}. When comparing this scenario to the case of infinite reset rate (Figure \ref{fig:3}, Panels B and C), we observe notable shifts in the positions of the Cusp points. For inhibitory coupling ($K = -2.0 < 0$), $\gamma$ increases from $0.2233$ to $0.3450$, while the corresponding value of $\eta_0$ decreases from $1.2886$ to $0.9381$.  In contrast, for excitatory coupling ($K = 2.0 > 0$), $\gamma$ decreases from $0.6347$ to $0.5841$, and the corresponding value of $\eta_0$ increases from $-3.6083$ to $-1.7305$.

	These results show that, even with relatively small reset rates, increasing the proportion of resetting neurons (i.e., increasing $\gamma$) enhances the stability of the system described by Eq. (\ref{4.18}). In the uniform rest state, an increase in resetting size expands the range of the reststate. However, during the uniform spiking state, under a smaller reset rate, the final averaged firing rate remains largely unaffected by changes in the proportion of resetting neurons.
	
	\subsubsection{Varying parameter $\lambda~\&~ \eta_0$}\label{subsubsec432}
	
	We now explore the following question: can a sufficiently large reset rate $\lambda$ approximate the behavior of the system in the infinite reset rate limit? 
	
\begin{figure}[htp]
		\centering
		\includegraphics[width=0.8\textwidth]{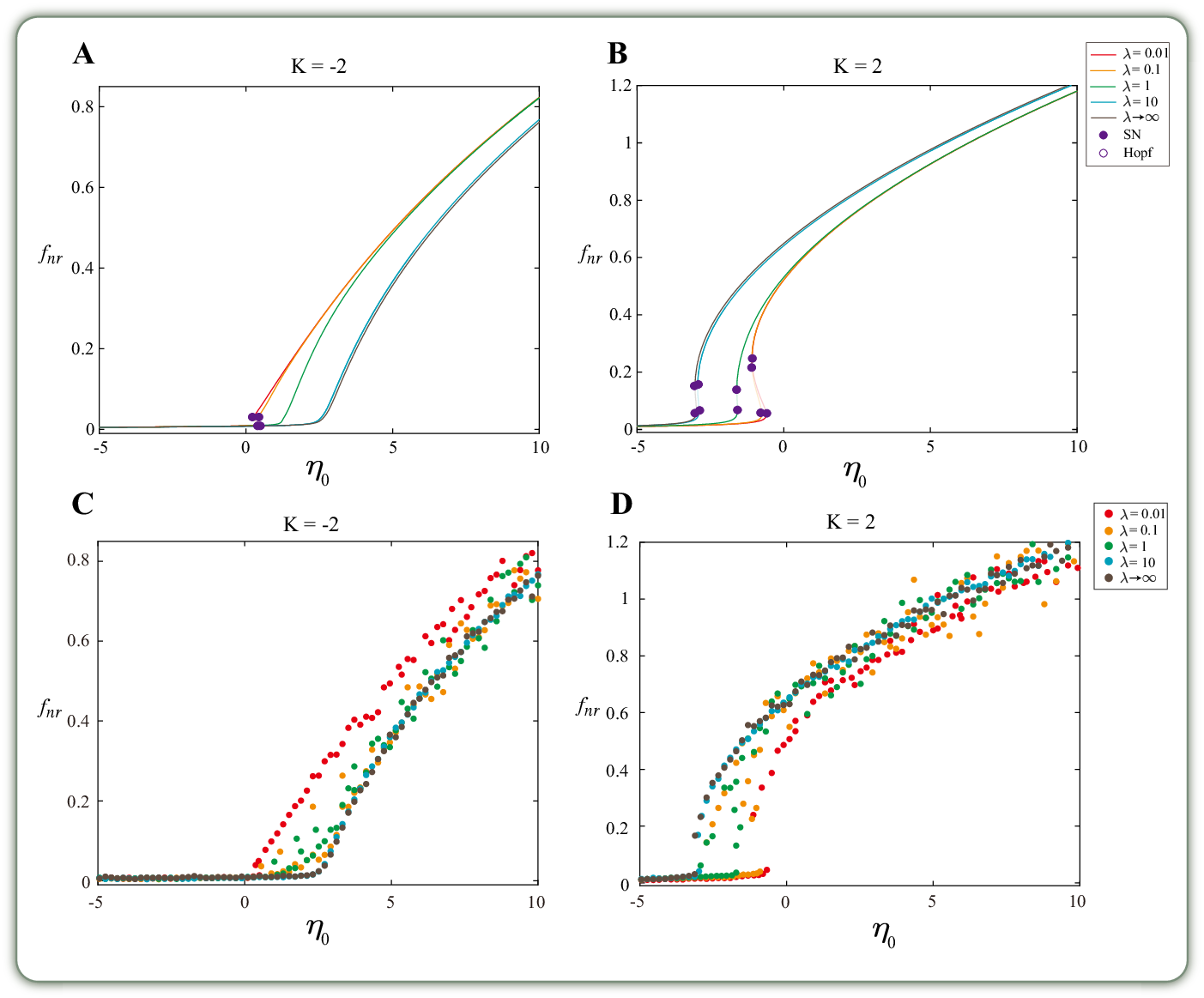}
		\caption{Plane A-B: Codimension-one bifurcation diagrams with $\eta_0$. $\gamma=0.5, \Delta=0.1$. For finite $\lambda$, the results are obtained by solving Eq. (\ref{4.18}), while for $\lambda\rightarrow\infty$, results are derived from solving Eq. (\ref{3.31}) simultaneously. Dark colors represent stable equilibria, while light colors represent unstable equilibria. Plane C-D Numerical simulation results for $\gamma = 0.5$, averaged $100$ realizations of $N = 1000$ and the integration time step equal to 0.01. See Table \ref{t3} for data on specific bifurcation points.}
		\label{fig:9}       
	\end{figure}

	To investigate this, we fix $\gamma=0.5$, meaning that half of the neurons are undergoing resetting. We then examine how varying the reset rate $\lambda$ affects the bifurcation behavior with respect to the parameter $\eta_0$ under both inhibitory and excitatory coupling, as illustrated in Figure \ref{fig:9}. 
	For comparison, we also plot the case where $\lambda \rightarrow \infty$ by solving Eq. (\ref{3.31}). When $\lambda$ is very small ($0.01$), the system's behavior closely follows the base evolution. However, as $\lambda$ increases to $1.0$, we observe that the lower branch of the bifurcation curve begins to shift. At higher firing rates, the system continues to adhere to the base evolutionary trend. When $\lambda$ reaches $10.0$, the entire bifurcation curve shifts significantly, no longer following the base evolution and instead approaching the behavior observed in the infinite reset rate case, as indicated by the brown curves in Figure \ref{fig:9}.
	
	The numerical simulation results show deviations from the theoretical results when $\lambda = 0.1, 1.0$. However, at $\lambda = 10.0 $, the results closely align with those expected in the asymptotic, infinite case. 
	
\begin{figure}[htp]
		\centering
		\includegraphics[width=0.8\textwidth]{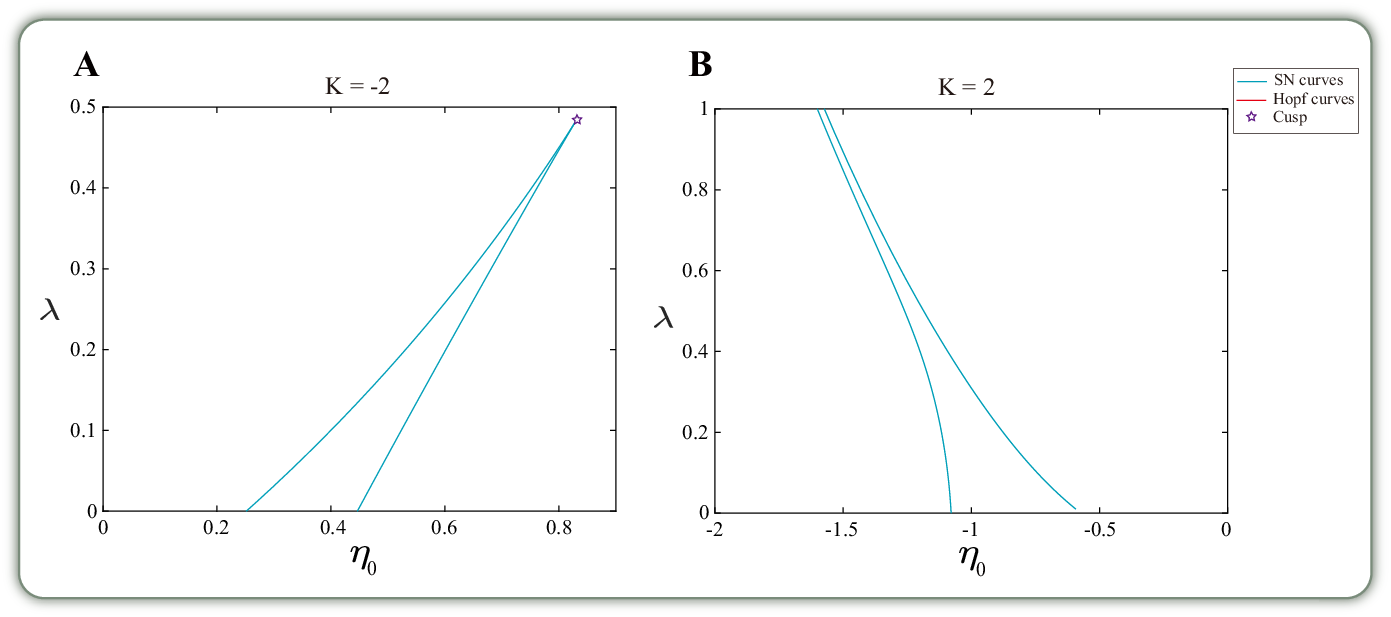}
		\caption{Codimension-two bifurcation in the plane $(\eta_0, \lambda)$. $\gamma=0.5, \Delta=0.1$. Legend is the same as Figure \ref{fig:3}. See Table \ref{t4} for data on specific bifurcation points.}
		\label{fig:10}       
	\end{figure}

	Then, we also plot the bifurcation curves in the parameter plane $(\eta_0, \lambda)$ . The results, shown in Figure \ref{fig:10}, highlight the changing system dynamics as $\lambda$ increases.

	In the case of inhibitory coupling, as $\lambda$ increases, the two saddle-node (SN) bifurcation curves eventually disappear at a Cusp point. However, for excitatory coupling, no such Cusp point is observed. Interestingly, even in the infinite reset rate case, the system does not achieve a fully stable equilibrium curve. As shown in Figures \ref{fig:3} C and \ref{fig:8} B, the Cusp points occur around $\gamma = 0.6$, indicating that approximately $60\%$ of the neurons must be reset to achieve a particular transition in system dynamics.
	
	In summary, we find that a reset rate of $\lambda=10.0$ already provides a close approximation to the infinite reset rate case. Therefore, in the simulation of full theta neuron network (Eq. \ref{2.1}), setting $\lambda$ as low as $10.0$  is sufficient to capture the essential dynamics of the infinite reset rate scenario. In other words, when $\lambda$ reaches $10.0$, the mean-field dynamics of the reset and non-reset subsystems almost entirely align with those of a single non-reset subsystem in the infinite reset rate case.
	
	\subsubsection{Varying parameter $\gamma~\&~ K$}\label{subsubsec433}
	
\begin{figure}[htp]
		\centering
		\includegraphics[width=0.8\textwidth]{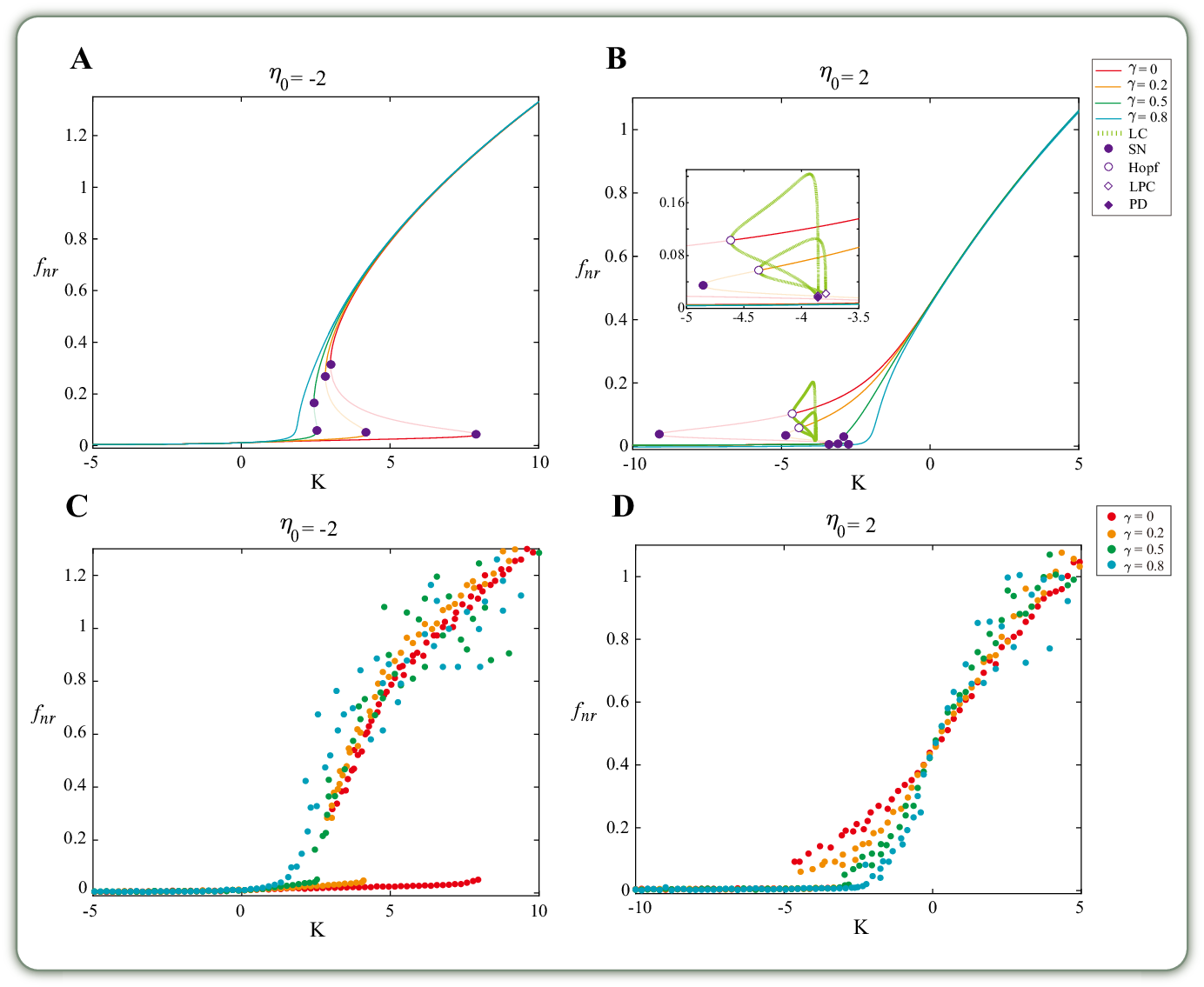}
		\caption{Plane A-B: Codimension-one bifurcation diagrams with $K$. $\lambda=1.0, \Delta=0.1$. Legend is the same as Figure \ref{fig:2}. Plane C-D: Simulations correspond to $N = 10000$ number of neurons with integration time step equal to $0.05$ and averaging over $100$ realizations. See Table \ref{t3} and \ref{t5} for data on specific bifurcation points.
		}
		\label{fig:11}       
	\end{figure}

	Next, we examine the influence of the global coupling strength $K$ as the bifurcation parameter, while keeping the reset rate $\lambda=1.0$. The corresponding results are shown in Figure \ref{fig:11}.
		
	Comparing Figure \ref{fig:5} B and C, we observe that for sufficiently large averaged firing rates (around $0.8$ for $\eta_0 = -2.0$ and $0.4$ for $\eta_0 = 2.0$), the equilibrium curves closely follow the base evolutionary trend. The transition from uniform rest to uniform spiking remains consistent with the results obtained for the infinite reset rate case. Specifically, for negative $\eta_0$ , the critical transition value decreases as $\gamma$  increases, whereas for positive $\eta_0$, the critical transition value increases.
	
	Notably, under infinite reset rate conditions, resetting more than half of the neurons ($\gamma = 0.5$) is sufficient to produce a fully stable equilibrium curve. However, under finite and relatively small reset rates, such as $\lambda=1.0$, approximately $80\%$ of the neurons must be reset to achieve the same level of stability.
	When more than half of the neurons are reset under a finite and relatively small reset rate ($\lambda =1.0$), the numerical results deviate from the theoretical predictions (green and blue curves). In all other cases, the numerical and theoretical results show good agreement.
	\begin{figure}[htp]
		\centering
		\includegraphics[width=0.8\textwidth]{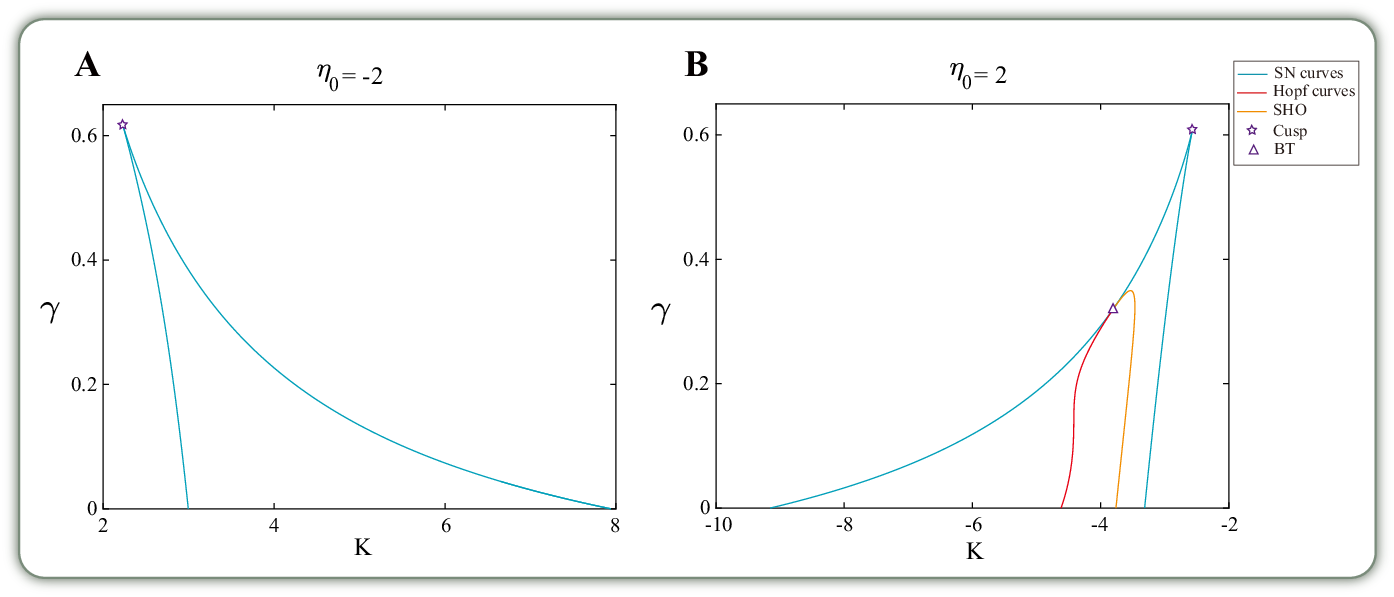}
		\caption{Codimension-two bifurcation in the plane $(K, \gamma)$. $\lambda=1.0, \Delta=0.1$. Legend is the same as Figure \ref{fig:3}. See Table \ref{t4} for data on specific bifurcation points.
		}
		\label{fig:12}       
	\end{figure}

	Additionally, we observe significant changes in the parameter plane $(K, \gamma)$, as illustrated in Figure \ref{fig:12}.
	Comparing the results with the infinite reset rate case in Figure \ref{fig:6}, we find that, for $\eta_0 = -2.0 < 0$, the Cusp point $\gamma$ increases from $0.4770$ to $0.6184$, while the corresponding $K$ increases from $1.3969$ to $2.2269$. For $\eta_0 = 2.0 > 0$, $\gamma$ increases from $0.3147$ to $ 0.6094$, while the corresponding $K$ decreases from $-2.2668$ to $-2.5635$. Furthermore, the BT points on the Hopf curve shift, with $\gamma$ increasing from $ 0.1385$ to $ 0.3220$, and $K$ increasing from $-3.8124$ to $-3.7894$.
	
	Thus, when using the global coupling strength $K$ as the bifurcation parameter, the averaged firing rate under a relatively small reset rate remains largely unaffected by increasing the number of reset neurons. However, the underlying transition mechanisms and the trend of the bifurcation curves still resemble those of the infinite reset rate case, although with slightly smaller magnitudes.
	
	\subsubsection{Varying parameter $\lambda~\&~ K$}\label{subsubsec434}
	
\begin{figure}[htp]
		\centering
		\includegraphics[width=0.8\textwidth]{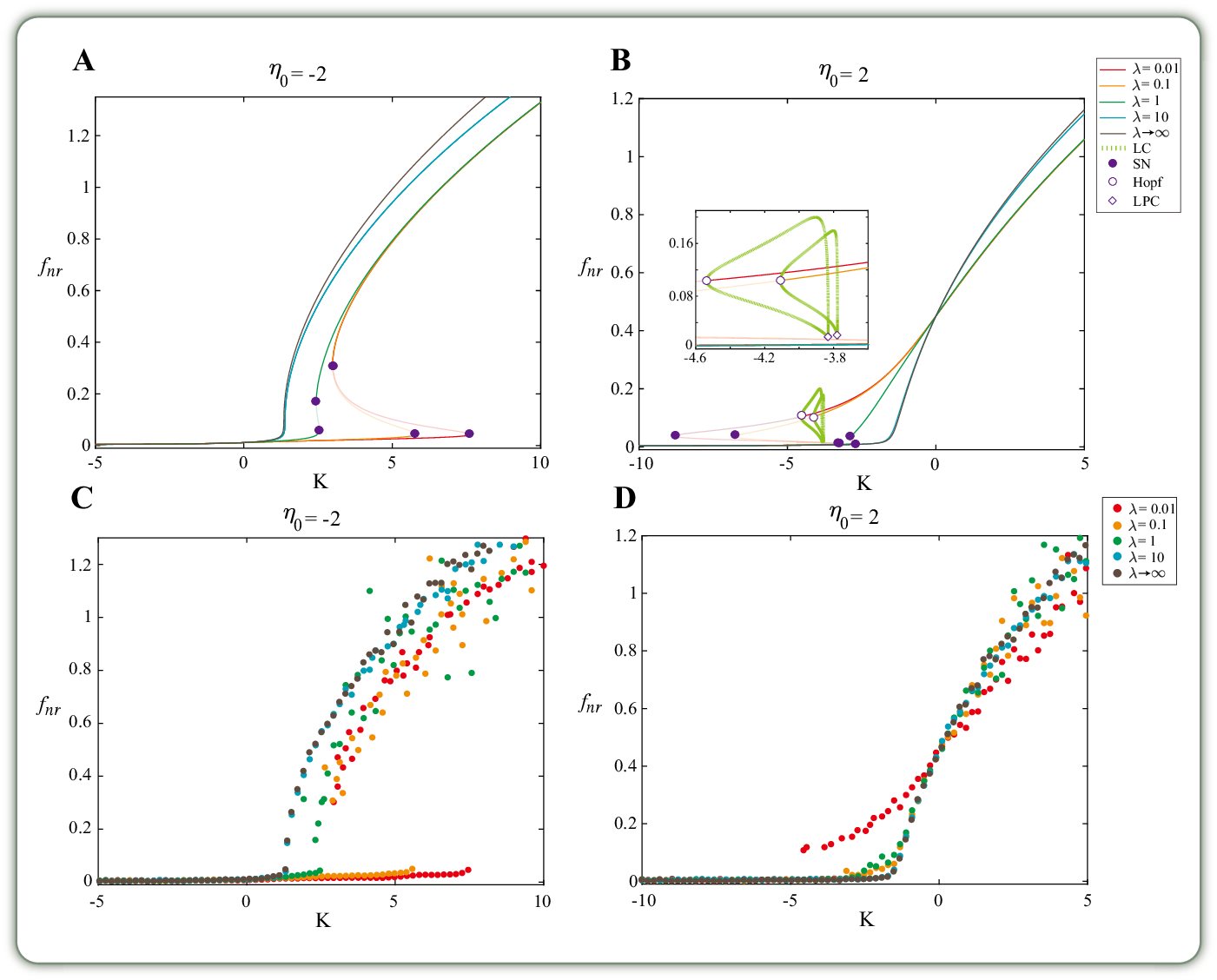}
		\caption{Plane A-B: Codimension-one bifurcation diagrams with $\eta_0$. $\gamma=0.5, \Delta=0.1$.  Legend is the same as Figure \ref{fig:9}. Plane C-D: Numerical simulation results for differnet finite $\lambda$, averaging over $100$ realizations of $10000$ neurons with the integration time step equal to $0.05$. See Table \ref{t3} and \ref{t5} for data on specific bifurcation points.
		}
		\label{fig:13}       
	\end{figure}

	Next, we investigate how the reset rate $\lambda$  influences the system’s bifurcation behavior, while keeping $\gamma=0.5$. The results are shown in Figure \ref{fig:13}.
	
	At very small reset rates, the system behavior aligns with the base evolutionary trend. As $\lambda$ is increased to $1.0$, the lower portion of the bifurcation curves begins to change. When $\lambda$ is further increased to $10.0$, the bifurcation curve deviates significantly from the base trend. By comparing these results to those obtained in the infinite reset rate case ($\lambda \rightarrow \infty$), we observe that as the reset rate increases, the resetting effect increasingly aligns with our expectations. However, for $\eta_0 = -2.0$, there is still some deviation in the upper right part of the bifurcation curve between $\lambda = 10.0$ and $\lambda\rightarrow\infty$. But, in the numerical simulations, $\lambda = 10.0$ also closely approximates the behavior expected in the infinite case.
	
\begin{figure}[!htp]
		\centering
		\includegraphics[width=0.8\textwidth]{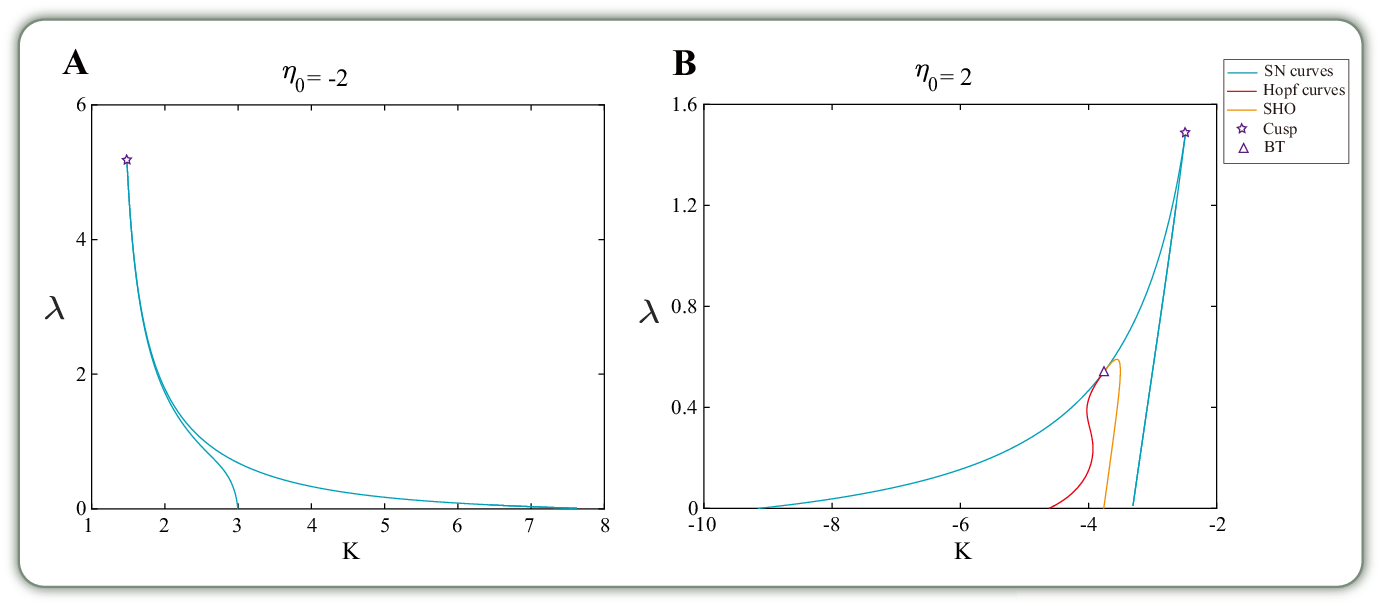}
		\caption{Codimension-two bifurcation in the plane $(K, \lambda)$. $\gamma=0.5, \Delta=0.1$. Legend is the same as Figure \ref{fig:9}. See Table \ref{t4} for data on specific bifurcation points.}
		\label{fig:14}       
	\end{figure}

	We also present the bifurcation curves in the parameter plane $(K, \lambda)$. As illustrated in Figure \ref{fig:14}. It is worth noting that when more than half of the neurons are reset, the system does not require a very large reset rate to stabilize the equilibrium curve fully.
	
	\subsubsection{Varying parameter $\gamma~\&~ \lambda$}\label{subsubsec435}
	
\begin{figure}[htp]
		\centering
		\includegraphics[width=0.8\textwidth]{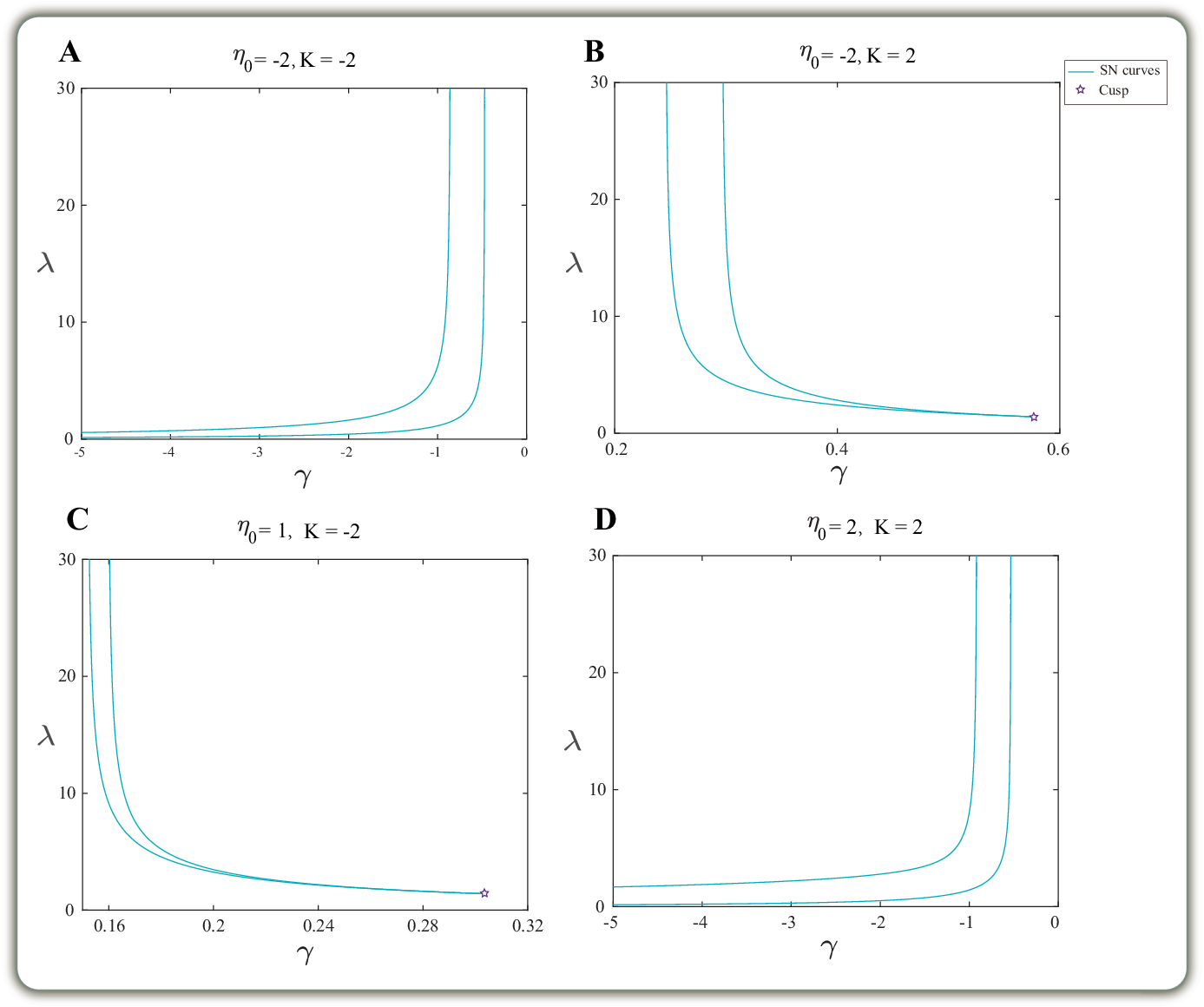}
		\caption{Codimension-two bifurcation in the plane $(\gamma, \lambda)$. $\Delta=0.1$. Legend is the same as Figure \ref{fig:9}. See Table \ref{t4} for data on specific bifurcation points.}
		\label{fig:15}       
	\end{figure}

	Building on the previous analysis, we now examine the behavior of the system within the parameter plane $(\gamma, \lambda)$, as illustrated in Figure \ref{fig:15}. 

	In this analysis, we consider four distinct scenarios: both median excitability ($\eta_0$) and global coupling strength ($K$) are positive, both are negative, and cases where one is positive while the other is negative. When $\eta_0$ and $K$ share the same sign (Planes A and D), we observe that the SN bifurcation curve unexpectedly extends into the negative range of the reset proportion, which contradicts our model assumptions. This suggests that configurations where $\eta_0$ and $K$ are both positive or negative may not adequately support the desired reset dynamics.
	
	When $\eta_0$ and $K$ have opposite signs, the behavior changes significantly. Here, the SN bifurcation curve gradually disappears at a Cusp point, indicating the collapse of the bistable region. This collapse facilitates a smooth and stable transition of the system from a uniform rest state to a uniform spiking state. These results imply that configurations with opposing values of $ \eta_0$ and $K$ better support an effective resetting mechanism, enhancing the system’s capacity for smooth state transitions. Such configurations may also improve stability and control within the reset dynamics of the theta neuron network, offering insights into optimized parameter setting.
	
	\section{Conclusion}\label{sec5}
	
	In this work, we investigated the macroscopic dynamics of theta neuron network subjected to partial resetting. Our focus shifts away from exploring how partial resetting influences global network synchronization. Instead, we analyze the stability of system equilibria using the average firing rate as an indicator, investigating the rich dynamics within the corresponding parameter space. The investigation was centered around two key scenarios: infinite reset rates and finite reset rates, with both exhibiting distinct effects on network dynamics, bifurcation structures, and stability.
	
	\subsection{Key findings}
	
	For the case of infinite reset rates, our analysis confirmed the theoretical predictions regarding the system's stability. As the reset rate approaches infinity, the reset subsystem becomes pinned to a fixed phase, exerting a continuous influence on the non-reset subsystem. The reset mechanism primarily facilitates a stable transition of the system from a uniform rest state to a uniform spiking state. This implies that the bistable region will disappear as the proportion of resetting neurons increases, due to the SN bifurcation curve ultimately colliding at a Cusp bifurcation at high $\gamma$ values. This version maintains clarity and emphasizes the effect of increasing the resetting neuron proportion. Although the system does not fully synchronize, it transitions into these uniform states, where neurons either maintain a steady resting potential or exhibit consistent spiking behavior. Our analysis showed that resetting enhances the stability of the averaged firing rate across different regions of parameter space. In particular, when the system is in a uniform spiking state, the neurons collectively exhibit a consistent spiking pattern that is strongly influenced by the resetting protocol, making the averaged firing rate more stable. In contrast, when the system is in a uniform rest state, the neurons remain largely inactive, and the reset mechanism ensures that this rest state is maintained, thereby stabilizing the averaged firing dynamics of the network. 
	
	When considering finite reset rates, the system's behavior becomes more complex. Finite resetting introduces stochasticity of the intermittent resect, which disrupts the OA ansatz in the infinite reset case. Thus, we analyze the system as two interacting populations: reset and non-reset subsystems. 
	With the reset rate $\lambda$ introduced as a control parameter, we performed a comprehensive bifurcation analysis across different regions of the parameter space. Our theoretical and simulation results indicate that a finite reset rate of $ \lambda = 10.0$ can produce effects similar to those of an infinite reset rate. However, when the reset rate is low, discrepancies emerge between theoretical predictions and numerical results. In particular, irregular jumps in the averaged firing rates are observed in the uniform spiking region, with substantial deviations from the predicted bifurcation curves. These deviations are especially pronounced when $50\%-80\%$ of neurons are reset, highlighting the limitations of the current analytical framework. Additionally, the reset mechanism is more effective when the median excitability and the global coupling strength have opposite signs, further supporting this configuration. 
	
	\subsection{Broader implications}
	This work provides several key insights that extend beyond the specific model of theta neuron network:
	
	1. \textbf{Resetting as a Stabilization Tool}: Both infinite and finite resetting offer a robust mechanism for stabilizing the averaged firing rate and controlling network dynamics. The ability of resetting to either enhance or suppress spiking behavior based on the coupling and excitability regimes makes it a versatile tool for modulating neural activity. This has implications for understanding biological neural processes where resets play a critical role, such as in memory, sleep, and pathological conditions like epilepsy.
	
	2. \textbf{Finite Reset Rates}: Our exploration of finite reset rates, an often overlooked aspect in neural network dynamics, shows that practical reset rates can approximate the behavior of infinite resetting. This finding is crucial for real-world systems, where infinite resetting is impractical, and offers insights into how biological systems might implement resets at finite rates while maintaining stability.
	
	3. \textbf{Applications in Neuromorphic Computing}: The findings from this study have potential applications in neuromorphic computing, where resetting mechanisms could be used to control the stability and function of artificial neural networks. Understanding how partial resetting influences spiking and resting states could lead to more efficient algorithms for managing network activity in hardware-based neural systems.
	
	While this study advances our understanding of resetting in neural networks, several open questions remain. The theoretical framework used for small reset rates could be further refined to capture the complex behavior observed in simulations. Additionally, the study focused on homogeneous networks, leaving room for exploration of more complex network topologies and heterogeneous neuron populations. Investigating the impact of stochasticity and noise on resetting behavior would also provide a more realistic model for biological neural systems.

	\section*{Acknowledgements}
	
	The authors extend their heartfelt thanks to the anonymous referee for their thorough review and insightful suggestions, which significantly improved the quality of this work. They are also grateful for the partial support provided by the National Natural Science Foundation of China (Grant No. 11872183) and the China Scholarship Council (Grant No. 202306150098).

	\appendix 
	\section{Details of numerical simulation of theta neuron network subject to subsystem resetting}\label{app1}
	
	\textbf{1. Parameter Setup}: Define key parameters, including the total number of neurons $N$, the fraction of resetting neurons $\gamma$, the reset rate $\lambda$, the global coupling strength $K$, median excitability $\eta_0$, and the total simulation time $T$. Choose $T$ based on $\lambda$ within a range from $20$ to $400$ to allow sufficient time for system dynamics to unfold. Set the degree of excitability variability to $\Delta = 0.1$, which controls the network’s heterogeneity.
	
	\textbf{2. Excitability Sampling}: Assign each neuron an intrinsic excitability $\eta_j$, sampled independently from a Lorentzian distribution $g(\eta)$. This step introduces natural variability across neurons, simulating the diverse excitability found in biological networks.
	
	\textbf{3. Initial Conditions}: Initialize the phase $\{ \theta_j(0) \}$ for each theta neuron, which also serves as the reset state for resetting neurons. Define the subset of resetting neurons as $b = \gamma N$. Following our setup, neurons labeled $j = 1, 2, \dots, b$ are classified as resetting neurons, while those labeled $j = b+1, b+2, \dots, N$ are non-resetting. For consistency in our implementation, all neurons are initially set to $\theta_j(0) = \pi, \forall j$. 
	
	\textbf{4. System Evolution}: From the initial conditions, allow the theta neuron network to evolve under the intrinsic dynamics specified by Eq. \ref{2.1}, for a random interval $\tau$sampled from an exponential distribution (Eq. \ref{2.6}). Numerically, this evolution is computed by integrating the system’s equations using a fourth-order Runge-Kutta algorithm. At the end of each interval $\tau$, the phases of the $b$ resetting neurons are instantaneously reset to $\pi$, while the phases of non-resetting neurons remain unchanged. This cycle of random-duration evolution followed by instantaneous resetting continues until the cumulative simulation time reaches $T$.
	
	\textbf{5. Multiple Realizations}: To account for variability in network dynamics due to disorder, repeat the simulation multiple times for each distinct excitability configuration $\{ \eta_j \}$. This allows reset intervals to vary across realizations. Averaging the outcomes across these realizations provides a stable estimate of network behavior under stochastic excitability variations.
	
	\textbf{6. Distribution Sampling}: To sample excitabilities $\eta_j$ from the Lorentzian distribution $g(\eta)$ use a uniform random-number generator, ensuring accurate representation of the specified distribution's properties.
	
	\section{Correlation data of the bifurcation points}\label{app2}
	
	In this section, we provide the supplementary data for all bifurcation points identified in our analysis. Please note $\text{SN}_l$ refers to the saddle-node bifurcation point on the lower branch of the bifurcation curve, while $\text{SN}_u$ corresponds to the point on the upper branch. Additionally, for special bifurcations such as Hopf, Cusp, and BT, LPC, PD bifurcation points, we provide the normal form coefficients, offering deeper insights into the local dynamics near these bifurcations. These coefficients are essential for characterizing the nature and stability of the bifurcations and are crucial for precise mathematical modeling.
	
	To ensure the highest level of precision, all bifurcation points are reported with four significant digits. This level of accuracy allows for a more detailed examination of the bifurcation ranges and minimizes potential numerical errors in subsequent analyses or simulations.
	
	\begin{table} [htp]
		\centering
		\caption{The data on the bifurcation points for synchronization results.}
		\footnotesize
		\label{t6}
		\begin{tabular}{lll}
			\hline\noalign{\smallskip}
			Figure \ref{fig:1} & $K=-2.0$ & $K=2.0$ \\
			\noalign{\smallskip}\hline\noalign{\smallskip}
			$\gamma=0$ & 
			
			$\text{SN}_{l}(0.2515, 0.8821)$, $\text{SN}_{u}(0.4464,0.9662)$ &
			
			$\text{SN}_{l}(-1.0789, 0.1287)$, $\text{SN}_{u}(-0.5730,0.7426)$ \\
			
			$\gamma=0.2$ & 
			
			$\text{SN}_{l}(1.1846, 0.9151)$, $\text{SN}_{u}(1.1914,0.9447)$ &
			
			$\text{SN}_{l}(-1.8257, 0.1999)$, $\text{SN}_{u}(-1.5375,0.7273)$ \\
			
			$\gamma=0.5$ & \quad & $\text{SN}_{l}(-3.0243, 0.3663)$, $\text{SN}_{u}(-2.9746,0.6765)$
			\\
			
			\noalign{\smallskip}\hline
		\end{tabular}
	\end{table}
	
	\begin{table} [htp]
		\centering
		\caption{The data on the Codimension-one bifurcation points for $\lambda\rightarrow\infty$.}
		\footnotesize
		\label{t1}
		\begin{tabular}{lllll}
			\hline\noalign{\smallskip}
			Figure \ref{fig:2} & $K=-10.0$ & $K=-2.0$ & $K=2.0$ & $K=10.0$ \\
			\noalign{\smallskip}\hline\noalign{\smallskip}
			\makecell[l]{$\gamma=0$\\ \\ \\ \\} & \makecell[l]{$\text{SN}_{l}(13.5445, 0.0066)$\\$\text{SN}_{u}(2.2011, 0.0376)$\\Hopf$(12.8792, 0.3897)$\\$l_1=-0.0227<0$ Sup} & 
			
			\makecell[l]{$\text{SN}_{l}(0.4464, 0.0133)$\\$\text{SN}_{u}(0.2515,0.0269)$\\ \\ \\} &
			
			\makecell[l]{$\text{SN}_{l}(-0.5730, 0.0516)$\\$\text{SN}_{u}(-1.0789, 0.2483)$\\ \\ \\} &
			
			\makecell[l]{$\text{SN}_{l}(-2.4879, 0.0423 )$\\$\text{SN}_{u}(-10.8527, 0.5942)$\\ \\ \\} \\

			\makecell[l]{$\gamma=0.2$\\ \\ \\ \\} & \makecell[l]{$\text{SN}_{l}(15.1012, 0.0071)$\\$\text{SN}_{u}(7.0605, 0.0370)$\\Hopf$(14.3328, 0.3182)$\\$l_1=0.0024>0$ Sub} & 
			
			\makecell[l]{$\text{SN}_{l}(1.1914, 0.0170)$\\$\text{SN}_{u}(1.1846,0.0217)$\\ \\ \\} &
			
			\makecell[l]{$\text{SN}_{l}(-1.5375, 0.0546)$\\$\text{SN}_{u}(-1.8257, 0.2144)$\\ \\ \\} &
			
			\makecell[l]{$\text{SN}_{l}(-7.3480, 0.0429)$\\$\text{SN}_{u}(-13.3909, 0.5335)$\\ \\ \\} \\

			\makecell[l]{$\gamma=0.5$\\ \\ \\ \\} & \makecell[l]{$\text{SN}_{l}(17.8727, 0.0085)$\\$\text{SN}_{u}(14.3445, 0.0351)$\\Hopf$(16.4746, 0.1548)$\\$l_1=0.2392>0$ Sub} & 
			
			\quad &
			
			\makecell[l]{$\text{SN}_{l}(-2.9746, 0.0652)$\\$\text{SN}_{u}(-3.0243, 0.1498)$\\ \\ \\} &
			
			\makecell[l]{$\text{SN}_{l}(-14.6350, 0.0446)$\\$\text{SN}_{u}(-17.5449, 0.4204)$\\ \\ \\} \\

			\makecell[l]{$\gamma=0.8$\\ \\} & \makecell[l]{$\text{SN}_{l}(21.7797, 0.0133)$\\$\text{SN}_{u}(21.5848, 0.0269)$} & 
			
			\quad &
			
			\quad &
			
			\makecell[l]{$\text{SN}_{l}(-21.9064, 0.0516)$\\$\text{SN}_{u}(-22.4122, 0.2483)$} \\

			\noalign{\smallskip}\hline\noalign{\smallskip}
			Figure \ref{fig:5} & $\eta_0=-10.0$ & $\eta_0=-2.0$ & $\eta_0=2.0$ & $\eta_0=10.0$ \\
			\noalign{\smallskip}\hline\noalign{\smallskip}
			\makecell[l]{$\gamma=0$\\ \\ \\ \\} & 
			\makecell[l]{$\text{SN}_{l}(41.8701, 0.0405)$\\$\text{SN}_{u}(9.4009,0.5769)$\\ \\ \\} &
			
			\makecell[l]{$\text{SN}_{l}(7.9381, 0.0429)$\\$\text{SN}_{u}(2.9956, 0.3171)$\\ \\ \\} &
			
			\makecell[l]{$\text{SN}_{l}(-3.3123, 0.0100)$\\$\text{SN}_{u}(-9.1507, 0.0374)$\\Hopf$(-4.6165, 0.1030)$\\$l_1=0.6213>0$ Sub} & 
			
			\makecell[l]{$\text{SN}_{l}(-8.1258, 0.0071)$\\$\text{SN}_{u}(-43.0830, 0.0394)$\\ Hopf$(-8.5206, 0.3383)$\\$l_1=-0.0061<0$ Sup} \\

			\makecell[l]{$\gamma=0.2$\\ \\ \\ \\} & 
			\makecell[l]{$\text{SN}_{l}(13.6703, 0.0421)$\\$\text{SN}_{u}(7.8644, 0.4733)$\\ \\ \\ } &
			
			\makecell[l]{$\text{SN}_{l}(2.6275, 0.0510)$\\$\text{SN}_{u}(2.1518, 0.2252)$\\ \\ \\ } &
			
			\makecell[l]{$\text{SN}_{l}(-2.6680, 0.0127)$\\$\text{SN}_{u}(-3.0482, 0.0296)$\\ \\ \\} &
			
			\makecell[l]{$\text{SN}_{l}(-7.3917, 0.0079)$\\$\text{SN}_{u}(-14.0671, 0.0379)$\\ Hopf$(-7.8947, 0.2420)$\\$l_1=0.0554>0$ Sub} \\

			\makecell[l]{$\gamma=0.5$\\ \\ \\ \\} & 
			\makecell[l]{$\text{SN}_{l}(6.8110, 0.0468)$\\$\text{SN}_{u}(5.9976, 0.3173)$\\ \\ \\} &
			
			\quad&
			
			\quad&
			
			\makecell[l]{$\text{SN}_{l}(-6.2091, 0.0103)$\\$\text{SN}_{u}(-7.0144, 0.0330)$\\ Hopf$(-7.0081, 0.0362)$\\$l_1=6.5824>0$ Sub} \\

			\makecell[l]{$\gamma=0.8$\\ \\ } & 
			\makecell[l]{$\text{SN}_{l}(4.5555, 0.0688)$\\$\text{SN}_{u}(4.5433, 0.1368)$} &
			
			\quad&
			
			\quad&
			
			\quad \\
			\noalign{\smallskip}\hline
		\end{tabular}
	\end{table}
	
	\begin{table} [htp]
		\centering
		\caption{The data on the Codimension-two bifurcation points for $\lambda\rightarrow\infty$.}
		\footnotesize
		\label{t2}
		\begin{tabular}{lllll}
			\hline\noalign{\smallskip}
			Figure \ref{fig:3} & $K=-10.0$ & $K=-2.0$ & $K=2.0$ & $K=10.0$ \\
			\noalign{\smallskip}\hline\noalign{\smallskip}
			\quad & \makecell[l]{Cusp$(22.6220, 0.8447)$\\ \qquad$c=-1.3966$ \\BT$(18.5003, 0.6716)$\\  \qquad $a=-2.9782$\\ \qquad$b=-1.6484$} & 
			
			\makecell[l]{Cusp$(1.2886, 0.2233)$\\ \qquad $c=-1.3966$ \\ \\ \\ \\ } &
			
			\makecell[l]{Cusp$(-3.6083, 0.6347)$\\ \qquad $c=-2.3829$\\ \\ \\ \\ } &
			
			\makecell[l]{Cusp$(-24.9416, 0.9269)$\\ \qquad $c=-2.3829$\\ \\ \\ \\ } \\

			\noalign{\smallskip}\hline\noalign{\smallskip} 
			Figure \ref{fig:6} & $\eta_0=-10.0$ & $\eta_0=-2.0$ & $\eta_0=2.0$ & $\eta_0=10.0$ \\
			\noalign{\smallskip}\hline\noalign{\smallskip} 
			\quad  & 
			\makecell[l]{Cusp$(4.3969, 0.8338)$\\ \qquad $c=-2.3829$ \\ \\ \\ \\} &
			
			\makecell[l]{Cusp$(1.3969, 0.4770)$\\ \qquad $c=-2.3829$\\ \\ \\ \\} &
			
			\makecell[l]{Cusp$(-2.2668, 0.3147)$\\ \qquad $c=-1.3966$ \\BT$(-3.8124, 0.1385)$\\  \qquad$a=-2.9782$\\ \qquad$b=-1.6484$} &
			
			\makecell[l]{Cusp$(-5.2668, 0.7051)$\\ \qquad $c=-1.3966$ \\BT$(-6.8124, 0.5179)$\\  \qquad $a=-2.9782$\\ \qquad$b=-1.6484$} 
			\\
			\noalign{\smallskip}\hline
		\end{tabular}
	\end{table}
	
	\begin{table} [htp]
		\centering
		\caption{The data on the Codimension-one bifurcation points for finite $\lambda$.}
		\footnotesize
		\label{t3}
		\begin{tabular}{lll}
			\hline\noalign{\smallskip}
			Figure \ref{fig:7} & $K=-2.0$ & $K=2.0$  \\
			\noalign{\smallskip}\hline\noalign{\smallskip}
			$\gamma=0$ & $\text{SN}_{l}(0.4464, 0.0133), \quad \text{SN}_{u}(0.2515, 0.0269)$ & 
			
			$\text{SN}_{l}(-0.5730, 0.0516), \quad \text{SN}_{u}(-1.0789, 0.2483)$ \\

			$\gamma=0.2$ & $\text{SN}_{l}(0.7133, 0.0149), \quad \text{SN}_{u}(0.6571, 0.0243)$ & 
			
			$\text{SN}_{l}(-0.9795, 0.0555), \quad \text{SN}_{u}(-1.2512, 0.2133)$ \\

			$\gamma=0.5$ & \quad &
			$\text{SN}_{l}(-1.5736, 0.0703), \quad \text{SN}_{u}(-1.6008, 0.1384)$ \\

			\noalign{\smallskip}\hline\noalign{\smallskip}
			Figure \ref{fig:9} & $K=-2.0$ & $K=2.0$  \\
			\noalign{\smallskip}\hline\noalign{\smallskip}
			$\lambda=0.01$ & $\text{SN}_{l}(0.4541, 0.0133), \quad \text{SN}_{u}(0.2672, 0.0267)$ & 
			
			$\text{SN}_{l}(-0.5927, 0.0517), \quad \text{SN}_{u}(-1.0798, 0.2479)$ \\

			$\lambda=0.1$ & $\text{SN}_{l}(0.5235, 0.0138), \quad \text{SN}_{u}(0.3999, 0.0250)$ & 
			
			$\text{SN}_{l}(-0.7434, 0.0543), \quad \text{SN}_{u}(-1.0925, 0.2432)$ \\

			$\lambda=1.0$ & \quad &
			$\text{SN}_{l}(-1.5736, 0.0703), \quad \text{SN}_{u}(-1.6008, 0.1384)$ \\

			$\lambda=10.0$ & \quad & 
			$\text{SN}_{l}(-2.8878, 0.0653), \quad \text{SN}_{u}(-2.9370, 0.1495)$ \\

			$\lambda\rightarrow\infty$ & \quad & 
			$\text{SN}_{l}(-2.9746, 0.0652), \quad \text{SN}_{u}(-3.0243, 0.1498)$ \\

			\noalign{\smallskip}\hline\noalign{\smallskip}
			Figure \ref{fig:11} & $\eta_0=-2.0$ & $\eta_0=2.0$  \\
			\noalign{\smallskip}\hline\noalign{\smallskip}
			\makecell[l]{$\gamma=0$\\ \\} & \makecell[l]{$\text{SN}_{l}(7.9381, 0.0429), \quad \text{SN}_{u}(2.9956, 0.3171)$\\ \\} & 
			
			\makecell[l]{$\text{SN}_{l}(-3.3123, 0.0100), \quad \text{SN}_{u}(-9.1507, 0.0374)$ \\Hopf$(-4.6165, 0.1030), l_1=0.5705>0$ Sub}\\

			\makecell[l]{$\gamma=0.2$\\ \\} & \makecell[l]{$\text{SN}_{l}(4.2430, 0.0477), \quad \text{SN}_{u}(2.8238, 0.2718)$\\ \\} & 
			
			\makecell[l]{$\text{SN}_{l}(-3.1074, 0.0109), \quad \text{SN}_{u}(-4.8576, 0.0348)$ \\Hopf$(-4.3731, 0.0573), l_1=1.2785>0$ Sub}\\

			$\gamma=0.5$ & $\text{SN}_{l}(2.5470, 0.0628), \quad \text{SN}_{u}(2.4393, 0.1660)$ & 
			
			$\text{SN}_{l}(-2.7411, 0.0136), \quad \text{SN}_{u}(-2.8939, 0.0271)$ \\
			
			\noalign{\smallskip}\hline\noalign{\smallskip}
			Figure \ref{fig:13} & $\eta_0=-2.0$ & $\eta_0=2.0$  \\
			\noalign{\smallskip}\hline\noalign{\smallskip}
			\makecell[l]{$\lambda=0.01$\\ \\} & \makecell[l]{$\text{SN}_{l}(7.6256, 0.0429), \quad \text{SN}_{u}(2.9949, 0.3168)$\\ \\} & 
			
			\makecell[l]{$\text{SN}_{l}(-3.3067, 0.0100), \quad \text{SN}_{u}(-8.8058, 0.0370)$ \\Hopf$(-4.5383, 0.1036), l_1=0.3072>0$ Sub}\\

			\makecell[l]{$\lambda=0.1$\\ \\} & \makecell[l]{$\text{SN}_{l}(5.8054, 0.0440), \quad \text{SN}_{u}(2.9854, 0.3137)$\\ \\} & 
			
			\makecell[l]{$\text{SN}_{l}(-3.2558, 0.0102), \quad \text{SN}_{u}(-6.7404, 0.0351)$ \\Hopf$(-4.1098, 0.1038), l_1=0.2564>0$ Sub}\\

			$\lambda=1.0$ & $\text{SN}_{l}(2.5470, 0.0628), \quad \text{SN}_{u}(2.4393, 0.1660)$ & 
			
			$\text{SN}_{l}(-2.7411, 0.0136), \quad \text{SN}_{u}(-2.8939, 0.0271)$\\
			
			\noalign{\smallskip}\hline
		\end{tabular}
	\end{table}
	
	\begin{table} [htp]
		\centering
		\caption{The data on the Codimension-two bifurcation points for finite $\lambda$.}
		\footnotesize
		\label{t4}
		\begin{tabular}{lll}
			\hline\noalign{\smallskip}
			\quad & $K=-2.0$ & $K=2.0$ \\
			\noalign{\smallskip}\hline\noalign{\smallskip}
			Figure \ref{fig:8} & Cusp$(0.9381, 0.3450)$, $c=-1.2840$ & 
			
			Cusp$(-1.7305, 0.5841)$, $c=-2.2244$ \\

			Figure \ref{fig:10} & Cusp$(0.8337, 0.4872)$, $c=-1.2231$ & 
			
			\quad \\
			
			\noalign{\smallskip}\hline\noalign{\smallskip} 
			\quad  &  $\eta_0=-2.0$ & $\eta_0=2.0$ \\
			\noalign{\smallskip}\hline\noalign{\smallskip} 
			\makecell[l]{Figure \ref{fig:12} \\ \\}  & 
			\makecell[l]{Cusp$(2.2269, 0.6184)$, $c=-2.2041$ \\ \\ } &
			
			\makecell[l]{Cusp$(-2.5635, 0.6094)$, $c=-1.2555$\\ BT$(-3.7894, 0.3220)$,  $a=-2.8675$, $b=-1.3377$} \\

			\makecell[l]{Figure \ref{fig:14}\\ \\}  & 
			\makecell[l]{Cusp$(1.4813, 5.1877)$, $c=-2.3867$ \\ \\ } &
			
			\makecell[l]{Cusp$(-2.4847, 1.4896)$, $c=-1.2974$\\ BT$(-3.7529, 0.5426)$,  $a=3.0701$, $b=1.2189$} \\
			\noalign{\smallskip}\hline
		\end{tabular}
	\end{table}
	
	\begin{table} [htp]
		\centering
		\caption{The data on the bifurcation points for limit cycles.}
		\footnotesize
		\label{t5}
		\begin{tabular}{llll}
			\hline\noalign{\smallskip}
			
			Figure \ref{fig:2} & $\gamma=0$ & $\gamma=0.2$ & $\gamma=0.5$ \\
			\noalign{\smallskip}\hline\noalign{\smallskip}
			\makecell[l]{$K=-10.0$\\\\\\} & \makecell[l]{$\text{LPC}(12.4366, 0.6510_{u}/0.0267_{l})$ \\ \\PD$(12.4610, 0.0100)$\\\qquad period $=7.7041$ \\} & 
			
			\makecell[l]{$\text{LPC}_1(14.3338, 0.3826_{u}/0.2549_{l})$ \\ $\text{LPC}_2(14.1769, 0.4985_{u}/0.0211_{l})$\\
				PD$(14.1845, 0.0106)$\\\qquad period $=8.6613$} & 
			
			\makecell[l]{$\text{LPC}(17.1612, 0.0116)$ \\\\ PD$(17.1612, 0.0122)$\\\qquad period $=11.1293$ \\}   \\
			
			\noalign{\smallskip}\hline\noalign{\smallskip} 
			Figure \ref{fig:5} & $\gamma=0$ & $\gamma=0.2$ & $\gamma=0.5$ \\
			\noalign{\smallskip}\hline\noalign{\smallskip}
			\makecell[l]{$\eta_{0}=2.0$\\\\\\} & \makecell[l]{$\text{LPC}(-3.8542, 0.0142)$ \\PD$(-3.8542, 0.0145)$\\\qquad period $=13.3795$ \\} & \quad & \quad \\
			
			\makecell[l]{$\eta_{0}=10.0$\\\\\\} &
			\makecell[l]{$\text{LPC}(-8.6536, 0.5533_{u}/0.0227_{l})$ \\PD$(-8.6470, 0.0104)$\\\qquad period $=8.2480$ \\} &
			
			\makecell[l]{$\text{LPC}(-7.7996, 0.4409_{u}/0.0580_{l})$ \\
				PD$(-7.8015, 0.0114)$\\\qquad period $=9.7525$} & 
			
			\makecell[l]{$\text{LPC}(-7.0011, 0.0291)$  \\\\ \\}   \\
			
			\noalign{\smallskip}\hline\noalign{\smallskip} 
			Figure \ref{fig:11} & $\gamma=0$ & $\gamma=0.2$ & \quad \\
			\noalign{\smallskip}\hline\noalign{\smallskip}
			\makecell[l]{$\eta_{0}=2.0$\\\\\\} &
			\makecell[l]{$\text{LPC}(-3.8542, 0.0145)$ \\PD$(-3.8542, 0.0144)$\\\qquad period $=13.9088$} &
			
			\makecell[l]{$\text{LPC}(-3.7905, 0.0181)$ \\\\\\} & 
			
			\quad     \\
			
			\noalign{\smallskip}\hline\noalign{\smallskip} 
			Figure \ref{fig:13} & $\lambda=0.01$ & $\lambda=0.1$ & \quad \\
			\noalign{\smallskip}\hline\noalign{\smallskip}
			$\eta_{0}=2.0$ &
			$\text{LPC}(-3.8313, 0.0142)$  &
			
			$\text{LPC}(-3.7774, 0.0148)$ & 
			
			\quad     \\
			
			\noalign{\smallskip}\hline
		\end{tabular}
	\end{table}

\newpage

\bibliographystyle{iet}{}
\bibliography{ied.bib}

\end{document}